\documentclass[12pt]{amsart}
\usepackage{enumerate,amssymb,amsfonts,latexsym,psfrag,psfig}
\usepackage{amscd,amsmath}
\usepackage[dvips]{graphicx}

\usepackage{fullpage}

\newtheorem{thm}{Theorem}[section]
\newtheorem{cor}[thm]{Corollary}
\newtheorem{lem}[thm]{Lemma}
\newtheorem{prop}[thm]{Proposition}
\theoremstyle{remark}
\newtheorem{rem}{Remark}[section]

\theoremstyle{definition}
\newtheorem{defn}{Definition}[section]

\numberwithin{equation}{section}
\numberwithin{figure}{section}

 \newcommand{\pichere}[2]
{\begin{center}\includegraphics[width=#1\textwidth]{#2}\end{center}}

\newcommand{\lab}[3]{\psfrag{#1}[#3]{$\scriptstyle{#2}$}}

\font\nt=cmr7

\def\note#1
{\marginpar
{\nt $\leftarrow$
\par
\hfuzz=20pt \hbadness=9000 \hyphenpenalty=-100 \exhyphenpenalty=-100
\pretolerance=-1 \tolerance=9999 \doublehyphendemerits=-100000
\finalhyphendemerits=-100000 \baselineskip=6pt
#1}\hfuzz=1pt}

\renewcommand{\Im}{\operatorname{Im}}

\newcommand{\di}{\partial}

\newcommand{\ra}{\rightarrow}

\def\ssk{\smallskip}
\def\msk{\medskip}

\def\nin{\noindent}

\def\ssm{\smallsetminus}

\newcommand{\diam}{\operatorname{diam}}
\newcommand{\dist}{\operatorname{dist}}

\newcommand{\inter}{\operatorname{int}}
\renewcommand{\mod}{\operatorname{mod}}

\newcommand{\tl}{\tilde}

\newcommand{\ctg}{\operatorname{ctg}}
\newcommand{\arcctg}{\operatorname{arcctg}}
\newcommand{\orb}{\operatorname{orb}}
\newcommand{\HD}{\operatorname{HD}}

\newcommand{\id}{\operatorname{id}}

\newcommand{\const}{\mathrm{const}}
\def\loc{{\mathrm{loc}}}

\newcommand{\eps}{{\varepsilon}}

\newcommand{\de}{{\delta}}
\newcommand{\la}{{\lambda}}
\newcommand{\La}{{\Lambda}}
\newcommand{\si}{{\sigma}}

\newcommand{\Om}{{\Omega}}
\newcommand{\om}{{\omega}}

\newcommand{\bare}{{\bar\eps}}

\newcommand{\AAA}{{\mathcal A}}

\newcommand{\CC}{{\mathcal C}}

\newcommand{\II}{{\mathcal I}}
\newcommand{\FF}{{\mathcal F}}
\newcommand{\GG}{{\mathcal G}}

\newcommand{\HH}{{\mathcal H}}

\newcommand{\OO}{{\mathcal O}}

\newcommand{\QQ}{{\mathcal Q}}

\newcommand{\UU}{{\mathcal U}}

\newcommand{\WW}{{\mathcal W}}

\newcommand{\C}{{\mathbb C}}

\newcommand{\D}{{\mathbb D}}

\newcommand{\N}{{\mathbb N}}

\newcommand{\R}{{\mathbb R}}

\newcommand{\Z}{{\mathbb Z}}

\def\Bf{{\mathbf{f}}}
\def\Bg{{\mathbf{g}}}
\def\BG{{\mathbf{G}}}

\def\Bv{{\mathbf{v}}}
\def\Bz{{\mathbf{z}}}
\def\Bx{{\mathbf{x}}}

\def\BF{{\mathbf{F}}}
\def\BS{{\mathbf{S}}}

\def\BPsi{{\boldsymbol{\Psi}}}
\def\BPhi{{\boldsymbol{\BPhi}}}
\def\B0{{\mathbf{0}}}

\def\BV{{\mathbf{V}}}
\def\BR{{\mathbf{R}}}
\def\BG{{\mathbf{G}}}

\newcommand{\Jac}{\operatorname{Jac}}

\newcommand{\spec}{\operatorname{spec}}

\catcode`\@=12

\def\Empty{}
\newcommand\oplabel[1]{
  \def\OpArg{#1} \ifx \OpArg\Empty {} \else
  	\label{#1}
  \fi}
		
\newcommand{\comm}[1]{}
\newcommand{\comment}[1]{}

\begin{document}

\bigskip\bigskip

\title[H\'enon renormalization]{Renormalization in the H\'enon family, I: \\
  universality but non-rigidity \\}
\author{A. de Carvalho, M. Lyubich, M. Martens}

\address{University of S\~ao Paolo}
\address {Stony Brook University and University of Toronto}
\address{University of Groningen}

\dedicatory{Dedicated to Mitchell Feigenbaum on the occasion
 of his 60th birthday}

\date{\today}

\begin{abstract}
  In this paper geometric properties of infinitely renormalizable
real H\'enon-like maps $F$ in $\R^2$ are studied. It is shown that the
appropriately defined renormalizations $R^n F$ converge
exponentially to the one-dimensional renormalization fixed
point. The convergence to one-dimensional systems is at a
super-exponential rate controlled by the average Jacobian and a
universal function $a(x)$. It is also shown that the attracting Cantor
set of such a map has Hausdorff dimension less than 1, but
contrary to the one-dimensional intuition,  it is not rigid,
does not lie on a smooth curve,  and generically has unbounded geometry.
\end{abstract}

\maketitle

\thispagestyle{empty}
\def\IMSmarkvadjust{0 pt}
\def\IMSmarkhadjust{0 pt}
\def\IMSmarkhpadding{0 pt}
\def\IMSpubltext{Published in modified form:}
\def\SBIMSMark#1#2#3{
 \font\SBF=cmss10 at 10 true pt
 \font\SBI=cmssi10 at 10 true pt
 \setbox0=\hbox{\SBF \hbox to \IMSmarkhpadding{\relax}
                Stony Brook IMS Preprint \##1}
 \setbox2=\hbox to \wd0{\hfil \SBI #2}
 \setbox4=\hbox to \wd0{\hfil \SBI #3}
 \setbox6=\hbox to \wd0{\hss
             \vbox{\hsize=\wd0 \parskip=0pt \baselineskip=10 true pt
                   \copy0 \break%
                   \copy2 \break%
                   \copy4 \break}}
 \dimen0=\ht6   \advance\dimen0 by \vsize \advance\dimen0 by 8 true pt
                \advance\dimen0 by -\pagetotal
	        \advance\dimen0 by \IMSmarkvadjust
 \dimen2=\hsize \advance\dimen2 by .25 true in
	        \advance\dimen2 by \IMSmarkhadjust

%
%
  \openin2=publishd.tex
  \ifeof2\setbox0=\hbox to 0pt{}
  \else 
     \setbox0=\hbox to 3.1 true in{
                \vbox to \ht6{\hsize=3 true in \parskip=0pt  \noindent  
                {\SBI \IMSpubltext}\hfil\break
                \input publishd.tex 
                \vfill}}
  \fi
  \closein2
  \ht0=0pt \dp0=0pt
 \ht6=0pt \dp6=0pt
 \setbox8=\vbox to \dimen0{\vfill \hbox to \dimen2{\copy0 \hss \copy6}}
 \ht8=0pt \dp8=0pt \wd8=0pt
 \copy8
 \message{*** Stony Brook IMS Preprint #1, #2. #3 ***}
}

\SBIMSMark{2005/07}{August 2005}{}

\setcounter{tocdepth}{1}
\tableofcontents

\section{Introduction}

  Since the universality discoveries, made in the mid-1970's by
Feigenbaum~\cite{F1,F2} and, independently, by Coullet and Tresser
\cite{CT,TC}, these fundamental phenomena have attracted a great deal
of attention from mathematicians, pure and applied, and physicists
(see~\cite{Cv} for a representative sample of theoretical and
experimental articles in early 1980's on the subject). However, a
rigorous study of these phenomena has been surprisingly difficult and technically
sophisticated and so far has only been thoroughly carried out in the
case of one-dimensional maps, on the interval or the circle, with one
critical point (see~\cite{FMP,L,Ma,McM,S,VSK,Y} and references therein).

Rigorous exploration of universality for dissipative two-dimensional
systems was begun in the article by Collet, Eckmann and Koch
\cite{CEK}. It is shown in this article that the one-dimensional
renormalization fixed point $f_*$ is also a hyperbolic fixed point for
nearby dissipative two-dimensional maps: this explained (at least, at the physical level)
parameter universality observed in families of such systems.  A subsequent
paper by Gambaudo, van Strien and Tresser~\cite{GST} demonstrates that,
similarly to the one-dimensional situation, infinitely renormalizable
two-dimensional maps which are close to $f_*$ have an attracting
Cantor set $\OO$ on which the map acts as the adding machine. However,
the geometry of these Cantor sets and global topology of the maps in
question have not yet received an adequate deal of attention.

\msk In this paper we begin a more systematic study of the
geometry of infinitely renormalizable dissipative two-dimensional
dynamical systems.\footnote{Here only period-doubling renormalization will be
considered and we will refer to it simply as ``renormalization.''} 
What we have discovered is that for these maps universality features
(some of which have specific two-dimensional nature)
can coexist with unbounded geometry and  lack of rigidity
(which make them quite different from the familiar one-dimensional
counterparts).

We consider a class $\HH$ of H\'enon-like maps of the form
$$
       F\colon (x, y) \mapsto (f(x)-\eps(x,y), x),
$$
where $f(x)$ is a unimodal map subject of certain regularity assumptions,
and $\eps$ is small.  
If $f$ is renormalizable
then the renormalization of $F$ is defined as $RF= H^{-1} \circ (F^2|_U)
\circ H$, where $U$ is a certain neighborhood of the ``critical value'' $v=(f(0),
0)$ and $H$ is an explicit {\it non-linear} change of variables (\S~\ref{Henren}).
\footnote{The set-up in this article is different from that of~\cite{CEK}: 
a different normalization of H\'enon-like maps is used, and
renormalization is done near the ``critical value'' rather than the
``critical point'' using, at least initially, a non-linear change of variable. 
We found the theory quite sensitive to specific choices such as this.}

It is shown that the degenerate map $F_*(x,y):= (f_*(x), x)$, where
$f_*$ is the fixed point of the one-dimensional renormalization
operator, is a hyperbolic fixed point for $R$ with a one-dimensional
unstable manifold (consisting of one-dimensional maps) and that the
renormalizations $R^n F$ of infinitely renormalizable maps converge at
a super-exponential rate toward the space of unimodal maps 
(Theorem ~\ref{convergence}
and ~\ref{fixed}). For any infinitely renormalizable map $F$ of class
$\HH$ there exists a hierarchical family of pieces $\{B^n_\si\}$,
$2^n$ on each level, organized by inclusion in the dyadic tree, such
that $$
\OO=\OO_F = \bigcap_n \bigcup_\si B^n_\si 
$$ 
is an attracting Cantor set on which $F$ acts as the adding
machine (Corollary ~\ref{adding machine}). 
This recasts the results of~\cite{CEK,GST} in our setting.

Furthermore, the diameters of the pieces $B^n_\si$ shrink at least exponentially with rate
$O(\la^{-n})$, where $\la=2.6\ldots$ is the universal scaling factor of 
one-dimensional renormalization (Lemma ~\ref{contracting}).
This implies that
$$
       \HD(\OO)<\log 2 / \log \la < 1,
$$
which makes it possible to control distortion of the renormalizations
(Lemma ~\ref{distortion}). 
Ultimately, this leads to the following
asymptotic formula for the renormalizations (Theorem ~\ref{universality}): 
$$ 
R^n F(x,y)
= (f_n(x) -\, b^{2^n}\, a(x)\, y\, (1+ O(\rho^n)), \ x\, ), 
$$ 
where $f_n\to f_*$ exponentially fast, 
$$ 
b=b_F= \exp \int_\OO \log \Jac F
\, d\mu, $$ is the {\it average Jacobian} of $F$ (here $\mu$ is the
unique invariant measure on $\OO$ and the Jacobian is the absolute value of the determinant), $\rho\in (0,1)$, and $a(x)$ is a
{\it universal} function. This is a new universality feature of
two-dimensional dynamics:
 as $f_*$ controls the zeroth order shape of the renormalizations, 
$a(x)$ gives the first order control.

\msk

On the other hand, we will  show in the second half of the paper
that there are some striking differences between the
one- and two-dimensional situations (\S~\ref{quadratic change of variable}  --  \S~\ref{unbounded geometry sec}).
For example, {\it the Cantor set $\OO$ is not rigid} (Theorem ~\ref{opthol}). 
Indeed, if the
average Jacobians of $F$ and $G$ are different, say $b_F < b_G $, then
a conjugacy $h\colon
\OO_F\ra \OO_G$ does not admit a smooth extension to $\R^2 $:
there is a definite upper bound
$$ \alpha\leq \frac{1}{2} \left(1+\frac{\log b_G}{\log b_F}\right) <1
$$ 
on the H\"older exponent of
$h$. Thus, in dimension two, universality and rigidity phenomena 
do not necessarily coexist. The above estimate on the H\"older exponent of the conjugation also applies to degenerate maps (i.e., one-dimensional) $F$ 
giving the upper bound $1/2$ on the H\"older exponent of $h$.

\begin{rem} 
One can compare this non-rigidity phenomenon with non-rigidity of circle maps.
In 1961 Arnold constructed an analytic diffeomorphism of the circle with  irrational 
rotation number whose conjugation with the corresponding rigid rotation is not 
absolutely continuous, see ~\cite{Ar}, ~\cite{H}. 
However, this phenomenon is quite different from the one discussed here as
it is related to the unbounded combinatorics (Liouville rotation number)
of the circle diffeomorphism in question.
\end{rem}

It was even more surprising to us that generically 
the Cantor set $\OO$ does not
have bounded geometry and so is not quasiconformally
equivalent to the standard Cantor set (Theorem ~\ref{unbdgeomth}).%
\footnote{In fact,  it seems to be quite a challenge to construct a {\it single example} of a H\'enon map of the class
we consider whose Cantor set would have bounded geometry (see Problem 5 in \S \ref{problems}).
It seems to go against the common intuition as one can find quite a few results in the  literature 
obtained under the assumption  of bounded geometry, compare \cite{CGM,Mo}}
Even worse, the Cantor sets
of generic infinitely renormalizable H\'enon-like maps have unbounded geometry  in some places, 
but in some other places they have a universal bounded geometry which is similar to their 
one-dimensional counterparts. 
(For instance, around the tip we always  recover the  universal scaling factor.) 
Moreover, the Cantor set $\OO$  cannot be embedded
into a smooth planar curve (Theorem ~\ref{curve}).

These properties, so different from their one-dimensional
counterparts, come from a {\it tilting} and {\it bending phenomenon}: near the
``tip'' of H\'enon-like maps renormalization boxes are not
rectangles but rather slightly tilted and  bent parallelograms. This tilt
significantly affects the $b$-scale geometry of $\OO$. Since the
Jacobian $b$ is replaced with $b^{2^n}$ under the $n$-fold
renormalization, the geometry gets affected at arbitrarily small
scales. These phenomena are explored in \S\ref{non-rigidity}, 
\S\ref{unbounded geometry sec} and \S\ref{holder}. 

The bent of the boxes forces us to use {\it non-affine} change of variables to make renormalizations
converge to a universal limit. However, we show in 
Theorem ~\ref{quad} 
that appropriate  quadratic changes of coordinates would be sufficient.
The renormalization limit obtained by this means would not correspond 
to the fixed point of  the usual renormalization  around the critical point, 
but rather to the one around the critical value.

In \S \ref{lineflds} we show that a non-degenerate H\'enon-like map in question
does not have continuous invariant line fields on the Cantor set $\OO$
(Corollary ~\ref{no line fields}). 
It implies that contrary to the ``rigidity intuition'',
the Cantor set $\OO$ does not lie on a smooth curve. 
It also implies that the ${\mathrm {SL}}(2,\R)$-cocycle  $z\mapsto DF(z)/ \sqrt{\Jac F(z)}$  is 
non-uniformly hyperbolic over the adding  machine 
$F: \OO\ra \OO$ (Theorem ~\ref{cocycleTh}).                                                    
It seems to be previously unknown whether such cocycles exist.

\ssk
On the positive side, as we show in the Theorem ~\ref{holgeo} , 
the Cantor set $\OO$ has 
{\it H\"older geometry} in an appropriate meaning of this term.

\comm{*****
The renormalization operator for  H\'enon-like maps differs 
from the standard unimodal period-doubling renormalization in the sense that 
it uses non-linear changes of variables.  
This is needed for the renormalizations to be H\'enon-like maps again.

A {\it posteriori}, we notice that, up to {\it bending} and {\it tilting}, 
the non-linear
changes of coordinates have a universal limit. Finally, 
one can use quadratic changes of coordinates
and see exponential convergence to a universal map.
This renormalization limit does not correspond 
to the fixed point of  the usual renormalization  around the critical point, 
but rather to the one around the critical value.  }

\msk In the forthcoming Part II, the global
topological structure of infinitely renormalizable H\'enon maps will
be discussed.

\msk
To conclude, it should be mentioned that intensive investigation of
stochastic attractors in the H\'enon family has been carried out
during the past two decades by Benedicks, Carleson, Viana, Young, and
others (see \cite{BC,BDV,WY}). This study has been concerned with
stochastic maps with positive entropy, which are very different from
the zero entropy maps studied here. We hope that, similarly to
what has happened in the one-dimensional theory, the renormalization
point of view will shed new light on stochastic phenomena as well.

\msk {\bf Acknowledgment.} We thank Jun Hu for sharing with us his
viewpoint on H\'enon renormalization, which is reflected in
\S~\ref{top def} of this paper. We thank Charles Tresser for
infinite renormalization discussions, and the referee for many useful comments. 
We also thank all the foundations that have supported us in the course of  this work:
the Guggenheim Memorial Foundation,
Clay Mathematics Institute, FAPESP, NSERC and NSF.

\section{General notation and terminology}
Let  $\N=\{1,2,\dots\}$, 
$\Z_+ = \N\cup \{0\}$, $I=[-1,1]\subset\R$, and 
\\ $\D_r= \{ z\in \C\colon \ |z|<r\}$.\\ 
A {\it rectangle} in $\R^2$ or $\C^2$ will mean a rectangle with
vertical and horizontal sides.

The letters $x$ and $y$ will be used
not only for real variables but also for their complexifications.  
The partial derivatives will be denoted by $\di_x$, $\di_y$, $\di^2_{xx}$, etc.  

For a smooth self-map $F$ of $\R^2$ or $\C^2$,
$\Jac F$ stands for $|\det DF|$. 

The coordinate projections in $\R^2$ or $\C^2$ will be denoted by $\pi_1$ and 
$\pi_2$. 
We let $\FF^h$ and $\FF^v$ be respectively
the foliations by horizontal and vertical real or complex lines in $\R^2$ or $\C^2$. 
A self-map of $\R^2$ or $\C^2$ is {\it horizontal} if it preserves the
horizontal foliation $\FF^h$.

A smooth map $f$ of an interval  is called {\it unimodal} if it has
a single critical point.  In what follows, we will assume that
{\it all the unimodal maps under consideration have a 
non-degenerate critical point and have negative Schwarzian derivative}, see \cite{dMvS}.


A self-map $H$ of $\R^2$ (from some family under consideration)
is said to have {\it bounded nonlinearity} 
if it may be represented as $A\circ\Phi$, where $A$
is affine and $\| \Phi - \id \|_{C^2} \leq K$, 
where $K$ is independent of the particular map is question.

The notation ``$\dist$'' will be used for different metrics in
different spaces, as long as there is no danger of confusion.  The
$\sup$-norm in the space $\AAA_\Om^c$  of bounded holomorphic functions on 
$\Om \subset \C^n$ is denoted by $\|\cdot\|_\Om$, or, if there is no danger of
ambiguity, simply by $\|\cdot\|$.  If $\Om$ is symmetric with respect to
the real subspace $\R^n$, $\AAA_\Om$ stands for the real slice of
$\AAA^c_\Om$ consisting of functions that are real on the real
subspace.


A set $X$ is called {\it invariant} under a map $f$ if $f(X)\subset X$.
$A\Subset B$ means that $A$ {\it is compactly contained in} $B$, i.e., 
the closure $\bar A$ is a compact  subset of $B$.  
Notation $Q_1\asymp Q_2$ means, as usual, that $C^{-1}\leq Q_1/Q_2\leq C$
for some constant $C>0$.

For reader's convenience, more special notations are collected in \S \ref{list}.

\section{H\'enon renormalization}

In this section, after briefly recalling the main definitions of
one-dimen\-sional renormalization, the class of H\'enon-like maps is
introduced and renormalization for such maps is defined. First a {\em
renormalizable map} is defined and this definition parallels the
one-dimensional definition: a certain topological disk is invariant
under the second iterate of the map. To define the {\em
renormalization} of the map, we consider the second iterate
restricted to the invariant disk and apply an appropriate non-linear change of
coordinates in order to obtain a H\'enon-like map of the same class.

\subsection{Renormalization of unimodal maps}\label{unimodal sec} 

A unimodal map $f\colon I\ra I$ with critical point $c\in I$ is called
{\it renormalizable} if there exists a closed interval $J\subset \inter I$
containing the critical point such that $J\cap f(J)= \emptyset$ and
$f^2(J)\subset J$. 
Then $f^2\colon J\ra J$ is a unimodal map.

We choose $J_c=[f^4(c),f^2(c)] $ to be the smallest interval as above,  
and call $f^2\colon J_c\ra J_c$ appropriately rescaled
(to bring $J_c$ back to the unit size)
the {\it renormalization} $R_c f$  of $f$. 
This is the classical {\it period-doubling} renormalization, 
and this is the only  renormalization type discussed in this paper.
However, we will also use the operator $R_v$ in the discussion of period 
doubling renormalization. It is defined as follows. Let 
$J_v=[f^3(c),f(c)] $ to be the smallest closed 
interval invariant under $f^2$ which 
contains the critical value $f(c)$,  
and call $f^2\colon J_v\ra J_v$ appropriately rescaled 
(to bring $J_v$ back to the unit size)
the {\it renormalization} $R_v f$  of $f$.
The operator $R_v$ renormalizes around the ``critical value '' and  
$R_c$ around the ``critical point''.

Let $r\in \Z_+\cup\{\om\}$ and let $\UU^r$ denote the space of
$C^r$-smooth unimodal maps $f\colon I\ra I$ such that:

\ssk \nin
 (a) the critical
point is mapped to $1$ and $1$ is mapped to $-1$ and 

\ssk\nin
(b) there is a
unique expanding fixed point $\alpha \in (-1,1)$ with negative
multiplier. 

The subspace of renormalizable maps is denoted by $\UU^r_0$, 
and the renormalization operators $R_c, R_v\colon\UU^r_0 \ra \UU^r$
assign to each map their renormalizations.%

For $r\ge 3$, the renormalization
operator $R_c$ has a unique fixed point $f_*\in \UU^\om_0$. 
It satisfies the functional equation $f_* = \la  f_*^2 (\la^{-1} x)$, 
where  $\la= 2.6\dots$ is the {\it universal scaling factor}. 
We let $\si= \la^{-1}$.

The fixed point $f_*$ is hyperbolic under the renormalization operator, 
with a  codimension-one  stable manifold $\WW^s(f_*)$  consisting of 
infinitely renormalizable maps. For details, see~\cite{L} and
references therein. The operator $R_v$ has also a unique fixed point $f^*$
(see Lemma 2.4 of \cite{BMT}).

\subsection{H\'enon-like maps}
Consider two intervals, $I^h$ and $I^v$, and let  $B=I^h\times I^v$.
A smooth  map $F\colon B\ra\R^2 $ is called {\it H\'enon-like}
if it maps vertical sections of  $B$ to horizontal arcs, while the
horizontal sections are mapped to parabola-like arcs 
(i.e., graphs of  unimodal functions over the $y$-axis).
%
Examples of H\'enon-like maps are given by small perturbations of
unimodal maps of the form
\begin{equation}\label{eq2}
    F(x,y) = (f(x) - \eps(x,y),\, x),
\end{equation}
where $f\colon I^h\ra I^h$ is unimodal and $\eps$ is small. Note that,
in this case, 
$$
         \Jac F = \left| \frac{\di\eps}{\di y} \right|. 
$$ 
If $\di\eps/\di y
\neq 0$ then the vertical sections are mapped diffeomorphically onto
horizontal arcs, so that $F$ is a diffeomorphism onto a ``thickening''
of the graph $\Gamma_f = \{(f(x), x)\}_{x\in I^h}$ (Figure 3.1).
\begin{figure}[htbp]
\begin{center}
\lab{g}{\Gamma_f}{l} \lab{B}{B}{b}
\pichere{0.6}{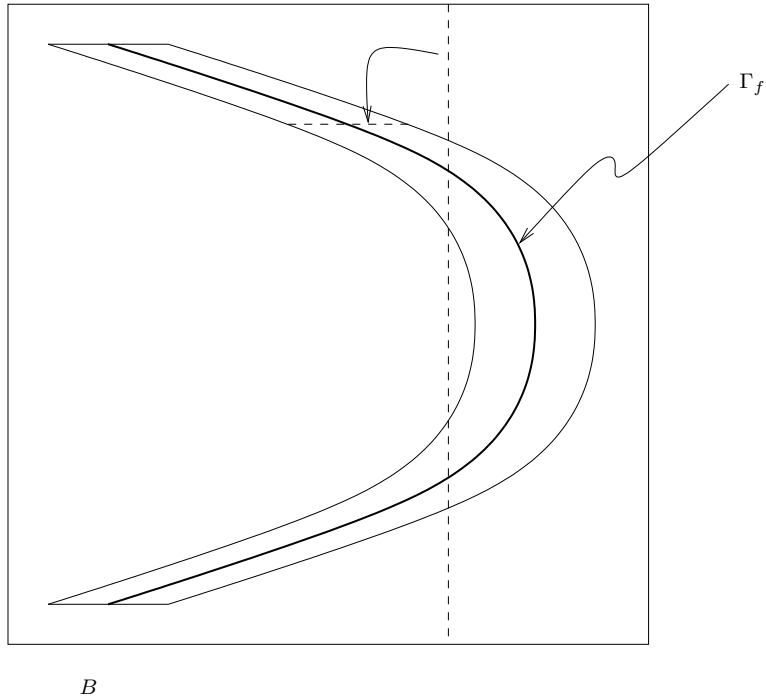}
\caption{A H\'enon-like map.}
\label{henon}
\end{center}
\end{figure}
In this case $F$ is a diffeomorphism onto its image which will be
briefly called a {\em H\'enon-like diffeomorphism}. 

The classical H\'enon family is obtained, up to affine
normalization, letting $f(x)$ be a quadratic polynomial and $\eps(x,y) = by$.

We will  use the abbreviation  $ F=(f-\eps,x) $
for equation (\ref{eq2}).  Thus, $F_f= (f, x)$ denotes the
degenerate H\'enon-like map collapsing $B$ onto $\Gamma_f$.

\footnote{ Usually, in particular in \S\ref{univ}, it is more convenient to 
consider unimodal maps only on 
their dynamical interval $[f^2(c), f(c)]$. 
However, without loss of generality we will assume that the unimodal maps 
have an extension defined on 
a {\it symmetric} interval bounded by an orientation preserving 
 fixed point and a preimage. We also assume, 
again without loss of generality, that all H\'enon-like maps have an extension 
containing the {\it regular} saddle point and its local stable manifold (compare 
\S\ref{top def} ).  }\label{exten}

\subsection{Spaces of maps}
Let $r\in \Z_+\cup \{\om\}$. The space of $C^r$-smooth H\'enon-like
maps $F\colon B\ra \R^2$ of the form (\ref{eq2}) is denoted by
$\HH^r$. Let $\UU^r$ be the space of unimodal maps as defined
above. In the real analytic case ($r=\om$), if $U\subset\C$ is a
neighborhood of $I$ and $\kappa>0$, then 
$\UU_{U,\kappa}\equiv\UU^\om_{U,\kappa}$
 denotes 
the subspace
of maps $f\in \UU_{U,\kappa}$ with critical point $c\in [-1,1-\kappa]$ 
which admit a holomorphic extension to
$U$ and  and can be factored as  $Q\circ \phi$, where $Q(x) =1 - x^2$ and
$\phi$ is an $\R$-symmetric  univalent map on $U$.                   
Since $\phi(c)=0$  and
$\phi(1) = \sqrt{2} $, this space of univalent maps is normal, so that
$\UU_{U,\kappa}$ is compact.\footnote{We fix once and for all a small $\kappa>0$ such that 
$c\in [-1,1-\kappa]$ for all maps of interest (like the renormalization fixed point and the infinitely renormalizable quadratic map), and we will suppress it 
from the notation.}

\ssk
Let $\Om^h, \Om^v\subset\D_2\subset\C$ be neighborhoods of $I^h,I^v$,
respectively, and let $\Om=\Om^h\times \Om^v\subset \C^2$.
 Let $\HH_\Om\equiv \HH^\om_\Om$ stand for the class of H\'enon-like maps
$F \in \HH^\om$ of form (\ref{eq2})
such that $f\in \UU_{\Om^h}$ and $\eps$ admits a holomorphic extension to $\Om$.  
The subspace of maps $F\in \HH_\Om$ with
$\|\eps\|_\Om\leq \bar\eps$ will be denoted by $\HH_\Om(\bar\eps)$. 
If $f$ in (\ref{eq2}) is fixed, we will also use the
notation $\HH_\Om(f, \bar\eps)$.

Realizing a unimodal map $f$ as a degenerate H\'enon-like map
$F_f$ yields an embedding of the space of unimodal maps
$\UU_{\Om^h}$  into the space of H\'enon-like maps $\HH_\Om$ making it
possible to think of $\UU_{\Om^h}$ as a subspace of $\HH_\Om$.

\subsection{Renormalizable H\'enon-like maps}
\label{top def}

\begin{figure}[htbp]
 \begin{center}
\lab{x0}{p_0}{l} \lab{x1}{p_1}{bl} \lab{x2}{p_2}{}
\lab{b0}{\beta_0}{lb} \lab{b1}{\beta_1}{}
 \pichere{0.6}{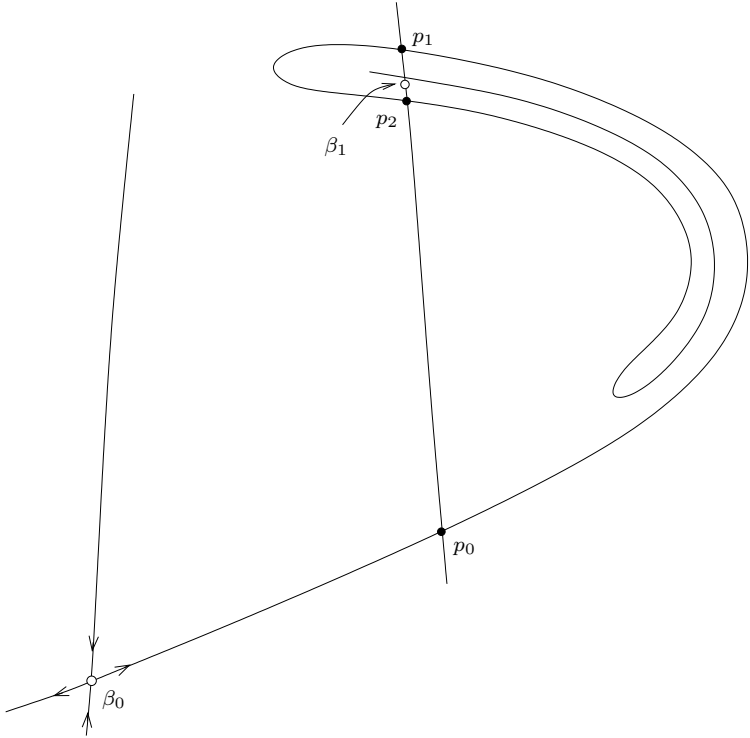}
 \caption{A renormalizable H\'enon-like map.}
 \label{invmflds}
 \end{center}
 \end{figure}

An orientation preserving H\'enon-like map is {\it renormalizable} if
it has two saddle fixed points --- a {\em regular} saddle $\beta_0$, with
positive eigenvalues, and a {\em flip} saddle $\beta_1$, with negative
eigenvalues --- such that the unstable manifold $W^u(\beta_0)$ intersects
the stable manifold $W^s(\beta_1)$ at a single orbit
(Figure~\ref{invmflds}).

For example, if $f$ is a renormalizable unimodal map with both fixed points repelling, 
then  a small H\'enon-like perturbation of type (\ref{eq2}) is a renormalizable
H\'enon-like map.

\msk Given a renormalizable map $F$, consider an
intersection point $p_0\in W^u(\beta_0)\cap W^s(\beta_1)$, and let
$p_n=F^n(p_0)$.  Let $D$ be the topological disk bounded by the arcs of
$W^s(\beta_1)$ and $W^u(\beta_0)$ with endpoints at~$p_0$ and~$p_1$.

\begin{lem}
The disk $D$ is invariant under $F^2$.
\end{lem}

\begin{proof} The boundary of $D$ consists of two arcs, $\ell^s\subset
W^s(\beta_1)$ and $\ell^u\subset W^u(\beta_0)$ both having $p_0$ and
$p_1$ for endpoints. Because $\beta_1$ is a flip saddle, $F^2(\ell^s)
\Subset \ell^s$ and there is a neighborhood $U\supset \ell^s$   with           \break 
 $F^2(U\cap D)\subset D$. If $F^2(D)$ were not contained
in $D $ then $F^2(\ell^u)$ would have to intersect the boundary $\ell^u\cup
\ell^s$ of $D$. The only possibility for this to happen would be that 
$F^2(\ell^u)$ intersects $\ell^s\setminus F^2(\ell^s)$. By hypothesis,
this intersection consists of points in the orbit of $p_0$. But this
would yield a contradiction, since $\ell^s\setminus F^2(\ell^s)$
contains only two points of the orbit of $p_0$, namely $p_0$ and
$p_1$, which are not in $F^2(\ell^u)$.
\end{proof}

\begin{defn}[Pre-renormalization]
The map $F^2|D$ is called a {\it pre-renormalization} of $F$.
\end{defn}

Assume now that $F$ is a small perturbation (\ref{eq2}) of a twice 
renormalizable 
unimodal map. In this case, there is a preferred intersection
point $p_0\in W^s(\beta_1) \cap W^u(\beta_0)$.  To define it, consider
the {\it local stable manifold} $W^s_\loc(\beta_1)$, the component of
the stable manifold $W^s(\beta_1)\cap B$ containing $\beta_1 $. If
$\eps$ is sufficiently small, then $W^s_\loc(\beta_1)$ is a nearly
vertical smooth arc. Let now $p_0$ be the {\it lowest} intersection
point of the unstable manifold $W^u(\beta_0)$ with
$W^s_\loc(\beta_1)$, so that the arc of $W^u(\beta_0)$ between
$\beta_0$ and $p_0$ does not intersect $W^s_\loc(\beta_1)$. This
determines the preferred pre-renormalization $F^2|D $ of $F$.

\subsection{The H\'enon renormalization operator}
\label{Henren}
We will now apply a carefully chosen  non-linear horizontal change of variables that
will turn the pre-renormalization into a H\'enon-like map of form (\ref{eq2}).

The pre-renormalization is not H\'enon-like, since it does
not map the vertical foliation to the horizontal one. 
However, it is not far from it:

\begin{lem}\label{complex slopes}
Let $f\in \UU_{\Om^h}$ with critical point $c$ and let $U \Subset
\Om^h\ssm \{c\}$ be an open set. There exist constants $C$ and
$\bar\eps > 0$, depending only on $\Om$ and $U$, such that
for any $F\in \HH_\Om(f,\bar\eps)$, the leaves of the foliation
$\GG=F^{-2} (\FF^h)$ in $U\times \Om^v$ are graphs over sub-domains of
$\Om^v$ with vertical slope bounded by $C \|\Jac F \|_\Om$.
\end{lem}

\begin{proof}
Since $\UU_{\Om^h}$ is a compact family of functions with a single critical point
$c\not\in \bar U$, we have $ \kappa := \min_{x\in U} |Df(x)| >0$, 
where $\kappa$ depends only on $\Om^h$.  
Letting
$r=\dist(\di U, \di \Om^h)$, if $\|\eps\|_\Om < \bare := \kappa r/2$,
then
\begin{equation}\label{x-der}
    \| \di\eps / \di x \|_{U\times \Om^v} <  \kappa/2.
\end{equation}

Since the foliation $F^{-2} (\FF^h)$ is given by the level
sets
$$
      f(x) - \eps(x,y) = \const
$$
it follows from the Implicit Function Theorem and (\ref{x-der}) that these level
sets are holomorphic graphs over sub-domains of $\Om^v$ with slopes satisfying
$$
   \left| \frac{\di x}{\di y} \right| =
     \left| \frac{\di\eps}{ \di y} \left( f'(x) - \frac{\di\eps}{\di x}
\right)^{-1} \right| \leq
      \frac{2}{\kappa} \left |\frac{\di\eps }{ \di y}\right| = \frac{2}{\kappa}
\Jac F (x,y).
$$ \end{proof}

For $U'\Subset U$, let $\Om' \subset \Om$ be the
{\it saturation} of $U'$ by  the leaves of the foliation $\GG\equiv F^{-2}(\FF^h)$, 
that is, $\Om'$ is the union of all leaves of $\GG$
 that intersect $U'$. 

\begin{cor}\label{U'}
If $U'\Subset U$ is an open set such that
$$
         \dist(\di U', \di U)> C \|\Jac F\|\diam \Om.
$$
then the leaves of $\GG$ that intersect $U'$ are holomorphic graphs
over $\Om^v$.
\end{cor}

Select neighborhoods $U'\Subset U\Subset \Om^h$ as above so that they
contain the interval $[\alpha,1]$ and $f|_U$ is an expanding
diffeomorphism with bounded non-linearity, with the bounds depending
only on $\Om$ and $U$. This is possible by compactness of $\UU_{\Om^h}$ and
because unimodal maps with negative Schwarzian derivative are expanding
on the interval $[\alpha,1]$.

\begin{lem}\label{eps2}
Given $U,U', \Om, \Om',  \GG$ as above, there exist $\bar\eps>0$, $C>0 $,
and a domain $V\ni c$ with the following  properties.  Consider a
H\'enon-like map $F=(f-\eps, x) \in \HH_\Om(f,\bar\eps)$ and define the
horizontal diffeomorphism
\begin{equation}\label{H}
  H(x,y) = (f(x) - \eps(x,y), y).
\end{equation}
Then there exists a unimodal map $g\in \UU_V$ such that $\| g -
f^2\|_V < C\bar\eps $ and $G:=H\circ F^2\circ H^{-1}$ is a
H\'enon-like map $(x,y)\mapsto (g(x)-\de(x,y), x)$ of class
$\HH_{V\times \Om^v}$ with $\|\de\|_{V\times \Om^v}\leq C\bar\eps^2 $.
\end{lem}

\begin{proof}
  Notice first that if $\eps$ is sufficiently small,
then all maps  $x\mapsto f(x)-\eps(x,y)$ are diffeomorphisms on
$U$ for any $y\in \Om^v$. Hence $H$ is a diffeomorphism  as well.

Let now
\begin{equation}\label{phi}
   \phi_y(x)= \phi(x,y) = f(x)-\eps(x,y),
\end{equation}
and
                $$ v(x)= -\eps(x, f^{-1} (x)).$$
A straightforward calculation gives us the following {\it Variational Formula}:
\begin{eqnarray}\label{VF}
  H\circ F^2\circ  H^{-1} (x,y) = \phi(\phi(x,\phi_y^{-1}(x)), x)  = 
\nonumber\\
  ( f^2(x)+ v(f(x)) + f'(f(x)) v(x) + O(\|\eps\|^2),\;  x), \label{VF2}
\end{eqnarray}
 which implies the assertion.
\end{proof}

\begin{rem}
 Note that $v$ is the restriction of the vector field $-\eps\, \di/\di x$
to the graph $\Gamma_f$, and $v\circ f + (f'\circ f) v$ is the first
variation of $f\mapsto f^2$ in the direction of $v$. Roughly speaking,
the two-dimensional variation of $f\mapsto f^2$ in the direction of
$-\eps$ coincides, to the first order, with its one-dimensional
variation in the direction of $v=-\eps|\Gamma_f$. In symbols: $
\de_{-\eps} (H\circ F_f^2\circ H^{-1}) = F_{\de_v f^2}.$
\end{rem}

\begin{rem}
 The residual term in (\ref{VF2}) involves second derivatives of
$\eps$, but in the holomorphic setting they are estimated by
$\|\eps\|$.
\end{rem}

\comm{*****
\begin{cor}
Take a neighborhood $U'\Subset U$ such that
$$
         \dist(\di U', \di U)> C \|\Jac F\|\diam \Om,
$$
and let $\Om' \subset \Om$ be the saturation of $U'$ with the
leaves of the foliation $\GG$. Then  there exists a holomorphic
horizontal diffeomorphism  $H\colon \Om' \ra
 U'\times\Om^v$ with $\|H - \id\|_{\Om'}= O(\|\Jac F\| ) $ that straighten
 the foliation $\GG $ to
the vertical foliation $\FF^v$.
\end{cor}

\begin{proof}
Notice first that since the vertical slope of the leaves of $\GG$
is bounded by $C\| \Jac F\|$, the leaves that intersect $U'$
never exit $\Om$. Hence they are
 holomorphic graphs over $\Om^v$.

Now,  let  $\gamma\colon \Om'\ra U'$ be the holonomy  along the
foliation $\GG$ which  slides a point $(x,y)\in \Om'$ along its
leaf of $\GG$ to a point
 $z\in U'$.
It is a holomorphic function such that  $ | \gamma(x,y) - x | \leq
C |\Jac F(x,y)|$. Then $H(x,y) = (\gamma(x,y), y)$ is a desired
map.
\end{proof}
*****}

\begin{defn}[Renormalization]
Let $J$ be the minimal interval such that $J\times I$ is invariant
under $G=H\circ  F^2\circ H^{-1}$, let $s \colon J\ra I$ be the
orientation-reversing affine rescaling, and let $\Lambda(x,y) = (s x, s y)$. 
Then the {\it renormalization}
 $RF$   is defined as $\La\circ G\circ \La^{-1}$ 
on the bidisk $\La ( V\times \Om^v)$.
\end{defn}

In the case of a degenerate map $F_f=(f,x)$ where $f$ is a
renormalizable unimodal map with critical point $c$, 
$J=[f^4(c),f^2(c)]$ is the same {\it dynamical interval} that we have used to
define the period doubling renormalization for unimodal maps.

Let us summarize the above analysis: 

\begin{thm}\label{R}
Given a domain $\Om\supset I$, there exist $\bar\eps>0$, $C>0$,
and  a neighborhood $s V$ of $I$ with the following properties.
Let $F=(f - \eps, x)$ be a renormalizable H\'enon-like map of
class $\HH_ \Om(\bar\eps)$. Then the renormalization $RF$ is a
H\'enon-like map of class $\HH_W(g,C\bar\eps^2)$,
where $W= \La(V\times \Om^v)$ and $g$ is a unimodal map such that
$\dist (R_c f,g) \leq C \bar\eps$. The change of variable $\La\circ
H$ conjugating $F^2$ (appropriately restricted) to $RF$ is an
expanding map with bounded non-linearity, with all bounds
depending only on $ \Om$ and $\bar\eps$. 
\end{thm}

\begin{rem}
Notice that if $F$ is close to the renormalization fixed point $F_*(x)
= (f_*(x), x)$ (see \S \ref{unimodal sec} and the next section), 
then the conjugacy $\La\circ H$ expands the
infinitesimal $l_\infty$-norm at least by factor $2.6$, 
as $\la=2.6\ldots$ is the dynamical scaling factor for the map $f_*$.
\end{rem}

\section{Hyperbolicity of the H\'enon renormalization operator}
\label{hyp}

In this section we show that the H\'enon renormalization operator
defined above has a hyperbolic fixed point 
\begin{equation}\label{F_*}
 F_*(x,y) = (f_*(x), x),
\end{equation}
 where $f_*$ is the fixed point of the one-dimensional renormalization
operator. We also show that, starting with an infinitely
renormalizable H\'enon-like map $F=(f-\eps,x)$ with $\eps$
sufficiently small, the renormalizations $R^n(F)$ converge
super-exponentially fast to the subspace of degenerate
(one-dimensional) maps and converge exponentially fast to the fixed
point $F_*$. It follows that the local unstable manifold $\WW^u(F_*)$
may be identified with the local unstable manifold $\WW^u(f_*)$, of the
one-dimensional renormalization operator, contained in the space of
unimodal maps, and that the local stable manifold $\WW^s(F_*)$ coincides
with the set of infinitely renormalizable H\'enon-like maps close to $F_*$. 

\bigskip

Let $\II_\Om(\bar\eps)$ and $\II_\Om(f,\bar\eps)$ denote the subspaces
of infinitely
renormalizable H\'enon-like maps (including degenerate ones) of
classes $\HH_\Om(\bar\eps)$ and $\HH_\Om(f, \bar\eps)$ respectively.

\begin{thm}\label{convergence}
Given a domain $\Om$, there is an $\bar\eps>0$ with the following
property: for $F \in \II_\Om(f,\bar\eps)$, there exists a domain
$V\subset \Om^h$ containing $I$ and a sequence of unimodal maps
$g_n\in \UU_V$ such that, for all $n\ge 0$,
$$ \| g_n - f_*\|_V \leq C \rho^n \|f-f_*\|_V$$ 
and 
$$\| R^n F - F_{g_n}\|_W =O(\bar\eps^{2^n}), $$ 
where $W=V\times \Om^v$ and $F_{g_n}=(g_n,x)$ is
the degenerate H\'enon-like map associated to $g_n$. All constants
depend only on $\Om$ and $\bar\eps$. The constant $\rho<1$ is
universal.
\end{thm}

\begin{proof}
 By the renormalization theory of unimodal maps,
it is possible to find  a domain $V \Subset \Om $ containing $I$ and  a
number  $N\in\N $ such that for any $N$ times renormalizable
unimodal map $f\in \UU_V$ the following holds:

\ssk\nin (i) $R_c^N f \in \UU_{V}$ and $ \dist
(R_c^N f,  f_*) <  (1/4) \dist ( f, f_*), $ where the distance is
associated with the norm $\|\cdot\|_V$.

It follows easily from the definition of the renormalization operator and compactness
of the space $\UU_V$  that

\ssk\nin (ii) There exists an $\bar\eps>0$ such that if $F\in
\II_{W}(f,\bar\eps)$ for some unimodal map $f\in \UU_{V}$, then
$f$ is $N$ times renormalizable.

\msk Take some  $\de> \dist( f, f_*)$. Let  $\bar\eps$ be so small
that property (ii) holds and $C\bar\eps< \min(1/2,\de/4)$,
where $C$ is the constant  from Theorem~\ref{R} applied to $R^N$. 
Let $g$ be a unimodal map approximating $R^N  F$ as given by Theorem~\ref{R}. Then
$$
  \dist(g, f_*) < \dist( g,R_c^N f) +  \dist(R_c^N f , f_*) <
$$
$$
   <  C\bar\eps +  (1/4)  \dist(f, f_*) < \de/2.
$$
Moreover,  $R^N F\in \HH_{W}(g, C\bar\eps^2)=\HH_{W}(g,
\bar\eps_1)$ with
$$ C\bar\eps_1 = (C \bar\eps)^2< (1/4) (\de/2).  $$ 
Hence it is possible to repeat the argument above with $R^N F$ in
place of $F$, $g$ in place of $f$, $\de/2$ in place of $\de$, and $
\eps_1$ in place of $\eps$.
In this way we construct inductively a sequence of
$N$-times renormalizable unimodal maps $g_k \in \UU_V$ such that
$\dist(g_k, f_*)<\de/2^k$ and $\dist(R^{Nk} F, g_k) = O(\bar\eps
^{2^k})$. The conclusion follows.
\end{proof}

By a standard trick (see, e.g., \cite[Prop. 3.3]{PS}), 
one can adapt the metric $\|\cdot\|$ to the
dynamics in such a way that $R$ becomes strongly contracting:

\begin{lem}\label{adapted metric}
  There is a metric on $\II_\Om(\bar\eps)$, equivalent to $\|\cdot\|_\Om$, and
$\rho\in (0,1)$ such that
$$
   \dist(RF, F_*) \leq \rho \dist (F, F_*)
$$
for all $F\in \II_\Om(\bar \eps)$.  \qed
\end{lem}

The space $\HH_\Om(\bar\eps)$ is naturally a real analytic Banach
manifold modeled on the space $\AAA_\Om$, with functions $\eps$
serving as local charts on $\HH_\Om(f, \bar \eps)$. It is obvious
from the definition that the renormalization operator $R\colon
\HH_\Om(\bar\eps)\ra \HH_\Om(\bar\eps)$ is real analytic.



By the unimodal renormalization theory,
the fixed point $f_*$ is a quadratic-like map on some domain $\Om_*\subset \C$
(see e.g, \cite{B} and references therein).
Moreover, $f_*$ is  a hyperbolic fixed point of $R_c$ in any space $\UU_V$
with  $V\Subset \Om_*$,

\begin{thm}\label{fixed}
Assume $\Om^h\Subset \Om_*$.  Then the map $F_*$ is the hyperbolic
fixed point for the H\'enon renormalization operator $R$ acting on
$\HH_\Om$, with one-dimensional unstable manifold $\WW^u(F_*)=\WW^u(f_*)$
contained in the space of unimodal maps. Moreover, the differential
$DR(F_*)$ has vanishing spectrum on the quotient $T \HH_\Om /
T\UU_{\Om^h}$.
\end{thm}

\begin{proof}
Let $E=T\HH_\Om/ T\UU_{\Om^h}$,  and let $A: E\ra E$ be the operator induced by  $DR(F_*)$. 
Then  Theorem~\ref{R} implies that $\|A^n\|=O(\bar\eps^{2^n})$, and hence 
$\spec(A)=\{0\}$.
\end{proof}

\begin{cor}
  The set  $\II_\Om(\bar\eps)$ of infinitely renormalizable H\'enon-like maps
coincides with the stable manifold
$$
  \WW^s(F_*) = \{F\in \HH_\Om(\bar\eps)\colon \ R^n F\to F_*\ \text{as} \ n\to
\infty \},
$$
which is a codimension-one real analytic submanifold in
$\HH_\Om(\bar\eps)$.   
\qed
\end{cor}

\begin{cor}
  For all $\Om$ and $\bar\eps$ as above,
  the intersection of $\II_\Om(\eps)$ with the H\'enon family
$$
  F_{a,b} \colon (x,y) \mapsto (f_a(x) - by, x)
$$
is a real analytic curve intersecting transversally the
one-dimensional slice $b=0$ at $a_*$, the parameter value for which $f_{a_*}$ is 
infinitely renormalizable.
\end{cor}

\begin{proof}
By the unimodal renormalization theory, the stable manifold
$\WW(f_*)= \WW^s(F_*)\cap\UU_\Om$ intersects transversally the quadratic family
$\QQ=\{f_a\}$ at a single point, $a_*$. 
By the hyperbolicity of the unimodal renormalization operator,
$R^n(\QQ)$ is close to $\WW^u(f_*)$ for big $n$'s.
Since $\WW^u(f_*)=\WW^u(F_*)$, the $R^n(\QQ)$ are transverse to $\WW^s(F_*)$ for big $n$'s as well.
It follows that $\QQ$, and hence the whole H\'enon family, is transverse to $\WW^s(F_*)$.  
\end{proof}

Let us finish this section with a complexification of the previous
results. Let $\HH^c_\Om(f_*, \bar\eps)$ stand for the space of maps of
form $F=(f_*-\eps,x)$, where $f_*\in \UU_{\Om^h}$ is the unimodal
renormalization fixed point and $\eps\in \AAA^c_\Om$ is a holomorphic
function on $\Om$ (not necessarily real on the real line) with
$\|\eps\|_\Om < \bar\eps$. This neighborhood of $F_*$ has a natural
complex structure inherited from $\AAA^c_\Om$, and the renormalization
operator $R$ extends to a holomorphic map on this space. 

\begin{thm}\label{complexification}
The degenerate map $F_*$ is a hyperbolic fixed point of the
renormalization operator $R$ acting on $\HH^c_
\Om(\bar\eps)$ with a codimension-one holomorphic stable
manifold $\II^c_\Om(\bar\eps) \equiv \WW^s_c(F_*)$, the
complexification of $\II_\Om(\bar\eps)=\WW^s(F_*)$.
\qed
\end{thm}

 The maps $F\in \II^c_\Om$ will still be called  infinitely renormalizable
(complex) H\'enon-like. Note that the renormalization of the
complex maps can be described geometrically in the same way as for
real maps, that is, as restriction of $F^2$ to an appropriate
bidisk, conjugating it by
 a horizontal map $H$
(given by the same formula) and rescaling.

\section{The critical Cantor set}
\label{Cantor}

Here we begin the study of the attracting set for infinitely
renormalizable H\'enon-like maps. As in dimension one, it is a Cantor
set on which the map acts like the dyadic adding machine. 
We show that its Hausdorff dimension is bounded from above by 0.73 and that it
depends holomorphically on the map. We will see in
Sections~\ref{non-rigidity} and \ref{unbounded geometry sec} 
that there are some fundamental differences
between these Cantor sets and their one-dimensional counterparts.

\bigskip

Consider an infinitely renormalizable complex H\'enon-like map $F\in
\II^c_ \Om(\bar\eps)$, where $\Om$ and $\bar\eps$ are selected so that
the previous results apply.

\subsection{Branches}\label{branches}

Let $\Psi^1_v\equiv \psi_v^1 := H^{-1}\circ \La^{-1}$ be the change of
variable conjugating the renormalization $RF$ to $F^2$ appropriately
restricted, and let $\Psi^1_c\equiv \psi_c^1= F\circ
\psi_v$. The subscripts $v$ and $c$ indicate that these maps are
associated to the critical {\it value} and the {\it critical} point,
respectively.

\begin{rem}\label{vertical fol preserved}
Note that  while the maps $\Psi_v^1$ preserve the horizontal foliation $\FF^h$,
the maps $\Psi^1_c$ preserve the vertical one, $\FF^v$.
Indeed, by definition (\ref{H}),
$H$ maps $F^{-1}(\FF^v)$ to $\FF^v$. Hence 
$$
   (\Psi_c^1)^{-1}(\FF^v)= \La\circ H(F^{-1}(\FF^v))=\FF^v.  
$$
\end{rem}

Similarly, let $\psi^2_v$ and $\psi^2_c$ be the corresponding
changes of variable for $RF$, let
$$
\Psi^2_{vv}= \psi^1_v\circ \psi^2_v, \quad \Psi^2_{cv}= \psi^1_c\circ
\psi^2_v, \quad \Psi^2_{vc}=\psi^1_v\circ\psi^2_c,\quad \dots.
$$
and, proceeding this way, construct, for any $n=1,2,\dots$, 
$2^n$ maps
$$ \Psi^n_w = \psi^1_{w_1}\circ\dots\circ \psi^n_{w_n}, \quad
w=(w_1, \dots, w_n)\in\{v,c\}^n. $$ 
The notation $\Psi^n_w(F)$ will also be used to emphasize dependence on the map $F$
under consideration, and we will let $W=\{v,c\}$ and
$W^n=\{v,c\}^n$ be the $n$-fold Cartesian product.

\begin{figure}[htbp]
 \begin{center}
\lab{F}{ F}{bl} 
\lab{RF}{RF}{bl} 
\lab{Rn-1F}{R^{n-1}F}{bl}
\lab{RnF}{R^nF}{bl} 
\lab{B1c}{B_c^1}{}
\lab{Q1c}{Q_c^1}{}
\lab{B1v}{B_v^1}{}
\lab{Q1v}{Q_v^1}{}
\lab{Qn-1c}{Q_c^{n-1}}{}
\lab{Qn-1v}{Q_v^{n-1}}{}
\lab{Qnc}{Q_c^n}{}\lab{Qnv}
{Q_v^n}{}
\lab{p1v}{\psi_v^1}{b}
\lab{p2v}{\psi_v^2}{b}
\lab{pn-1v}{\psi_v^{n-1}}{b}
\lab{pnv}{\psi_v^n}{b}
\pichere{1.0}{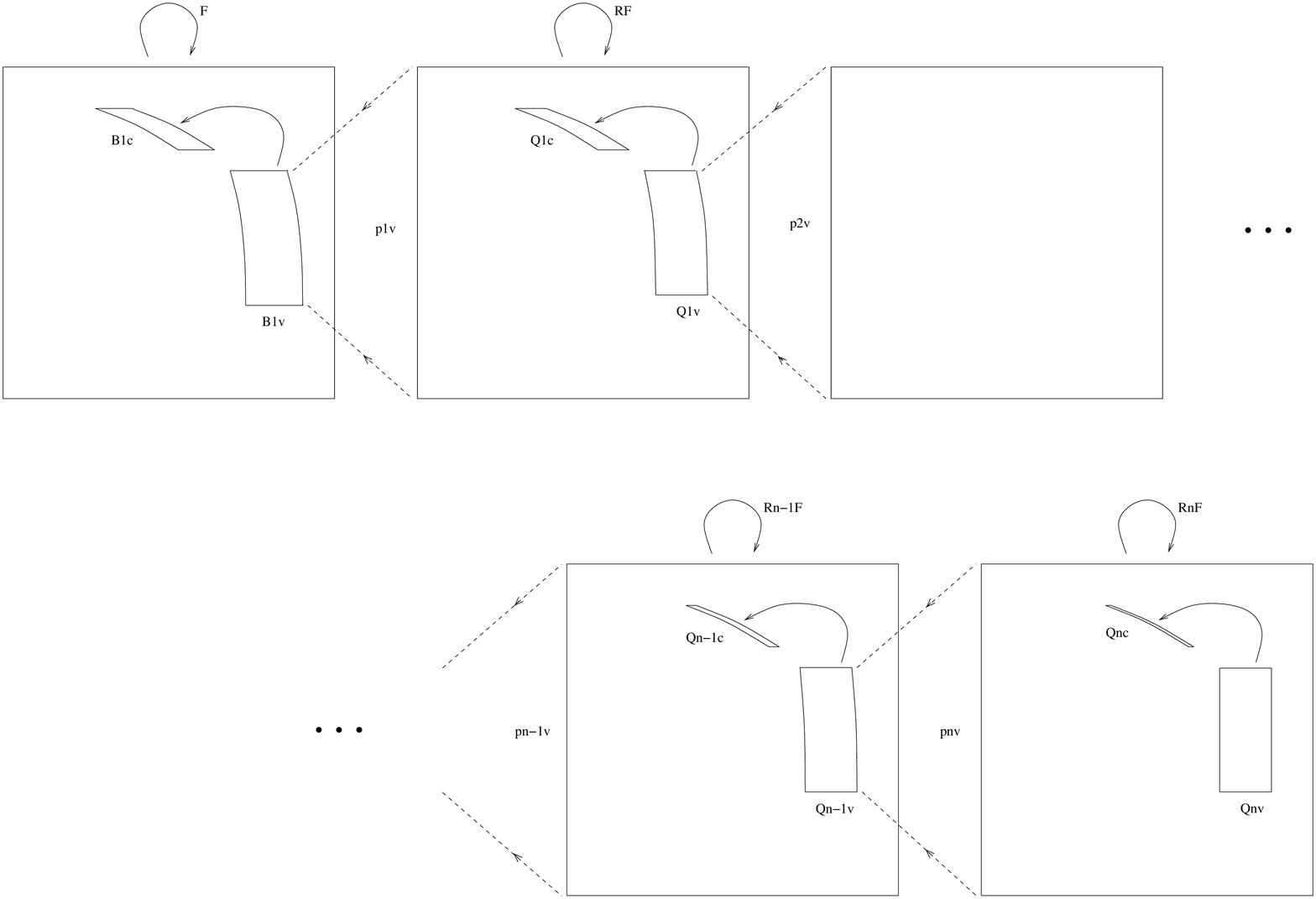}
\caption{The renormalization microscope}
\label{ren}
\end{center}
\end{figure}

Recall that $\si=\la^{-1}$ where $\la$ is the universal scaling factor.

\begin{lem}\label{contracting}
Let $F\in \II^c_\Om(\bar\eps)$,
$n\ge 1$, and $w\in W^n$. There exist $C>0$ and a domain 
in $\C^2$,
depending
only on $\Om$ and $\bar\eps$, on which the  holomorphic map $\Psi^n_w$ is 
defined and
$\| D\Psi^n_w\|\leq  C \sigma^n$.
\end{lem}

\begin{proof}
In the notation from equation (\ref{phi}) we have:  
$$ 
     H^{-1} (x,y) = (\phi_y^{-1} (x), y) \quad {\rm and} \quad   F\circ H^{-1} (x,y) = (x, \phi_y^{-1} (x)).
$$
The map $\phi_y^{-1}$ is
uniformly contracting on a neighborhood of the interval $J$, so that
$|\di\phi_y^{-1}/\di x|$ is bounded away from 1. On the other hand,
$\di \phi_y^{-1}/\di y $ is comparable with $\di\eps / \di y$, which
is small. It follows that the maps  $\psi_v =  H^{-1}\circ \La^{-1}$ and 
$\psi_c= F\circ H^{-1}\circ \La^{-1}$
 uniformly contracts the infinitesimal $l_\infty$-metric 
at least as strongly as $\Lambda^{-1}$, 
that is, by  a factor  $\si (1+O(\dist(F, F_*))$.

Since $R^nF \to F_*$ exponentially fast, the maps $\psi_{w_k}^k$, $w_k\in W$,
contract the infinitesimal $l_\infty$ normal by a factor $\si(1+O(\rho^k))$,
where $\rho\in (0,1)$. 
Hence the compositions $\Psi^n_w$ of these maps are
contracting by a factor $O(\si^n)$. 
\end{proof}

\subsection{Pieces}\label{pieces}

 Let us define $B_v^1\equiv B_v^1(F)=\psi^1_v(B)$
and $B_c^1\equiv B_c^1(F) = F(B_v^1)$. 
 Then $F(B_c^1)\subset B_v^1$. 
We will let $Q^n_w = B^1_w(R^n F)$, $n\in \Z_+$,
$w\in W$. Let $Q^\infty_w$ stand for the corresponding
pieces for the  degenerate limit  map (\ref{F_*}). Note that the
pieces $Q^n_w$  depend on $F$ while the pieces $Q^\infty_w$ do
not,  and that the piece $Q^\infty_c$ is in fact an arc on the parabola-like curve $x=f_*(y)$.

\begin{lem}\label{Q}
Let $F\in I^c_\Om(\bar\eps)$.  The pieces $Q^n_v$ and $Q^n_c$ have
disjoint projections to both of the coordinate axes. Moreover, they
converge exponentially, in the Hausdorff topology, to the pieces
$Q^\infty_v$ and $Q^\infty_c$, respectively.
\end{lem}

\begin{proof} The first statement follows easily from the definition
of renormalization. The second one follows from the exponential
convergence $R^n F \to F_*$.
\end{proof}

The sets $B^n_{w}\equiv B^n_{w}(F) = \Psi^{n}_w (B)$, where
$w\in W^{n}$, will be called {\em pieces}. They are closed
topological disks. For each $n\in \N$, there are $2^n$ such pieces and
forming the $n^{\text{th}}$-{\em generation} or 
$n^{\text{th}}$-{\em level} pieces.
$W^n$ can be viewed as the additive group of residues $\mod 2^n$ by letting 
$$ w
\mapsto\sum_{k=0}^{n-1} w_{k+1} 2^k,
 $$
where the symbols $v$ and $c$ are
interpreted as 0 and 1 respectively. Let $p\colon W^n\ra W^n$ be
the operation of adding 1 in this group.

\begin{lem}\label{permute}
\begin{enumerate}
\item[{\rm (1)}] The above families of pieces are nested:
$$
   B^n_{w\nu}\subset B^{n-1}_w, \quad w\in W^{n-1}, \ \nu\in W.
$$

\item[{\rm (2)}] The pieces $B^n_w$, $w\in W^n$, are pairwise
disjoint.

\item[{\rm (3)}] Under $F$, the pieces are permuted as follows.
$F(B^n_w) = B^n_{p(w)}$ unless $p(w) = v^n$. If
$p(w)=v^n$, then $F(B^n_w) \subset B^n_{v^n}$.
\end{enumerate}

\comm{ \nin (iv) In the latter case, the  piece $B^n_{v^n}$ is
contained in a $\|\eps\|^{2^n}$-neighborhood  of $F(B^n_\si)$. }

\end{lem}

\begin{proof} The first assertion holds by construction:
$$
B^n_{w\nu}= \Psi^{n}_{w \nu} (B)= \Psi^{n-1}_{w}\circ
\psi^{n}_{ \nu} (B)\subset B^{n-1}_{w}.
$$
The second follows by induction.
 For all maps under consideration we have by Lemma \ref{Q}
 that $B^1_v(F)$ and $B^1_c(F)$ are disjoint. Assume that the pieces of the
 $n^{th} $ generation
are disjoint for all maps under consideration. This implies that
the pieces $B^{n+1}_{w v}\subset B^1_v$, $w\in W^n$,  are pairwise disjoint,
as they are images of the disjoint pieces $B^{n}_{w }(RF)$  by the map $\psi_v^1$. 
Applying $F$, we see that the pieces $B^{n+1}_{wc}\subset B^1_c$, $w\in W^n$,
 are pairwise disjoint as well.
The assertion follows because  $B_c^1$ and $B_v^1$ are also disjoint.

\ssk Let us inductively check the third assertion. For $n=1$, we
have:
$$
 \mbox{ $B^1_c= F (B^1_v)$ and $F(B^1_c) = F^2(B^1_v)\subset B^1_v$.}
$$

Consider now  the pieces $B^n_w(RF)$, $w\in W^n$, of level $n$
for $RF$. Assume inductively that they are permuted by $RF$ as
required. Then the pieces $B^{n+1}_{v w} = \psi^1_v (B^n_w(RF))$, $w\in W^n$,
 are permuted in the same fashion under $F^2$.
Moreover, $B^{n+1}_{cw} = \psi^1_c (B^n_w(RF)) = F(B^{n+1}_{vw})$, 
and the conclusion follows.
\comm{*****

\ssk For the last assertion, let us consider the renormalization
$G_n \equiv R^n F $. Let $\eps_n$ be the distance from $G_n$  to
the corresponding unimodal map. Then the piece $Q^n_v$ is
contained in the $\eps_n$-neighborhood of $G_n(Q^n_v)$. Applying
the contracting map $\Psi^n_v$, we see that $B^n_{v^n}$ is
contained in the $\rho^n \eps_n$-neighborhood of
$F^{2^n}(B^n_{v^n})$. Since $\eps_n = O(\|\eps\|^{2^n})$ by Lemma
\ref{eps2}, we are done. 
*****}
\end{proof}

Furthermore, Lemma~\ref{contracting} implies:

\begin{lem}\label{boxes shrink}
There exists $C>0$, depending only on $\Om$ and $\bar\eps$, such that
for all $w\in W^n$, $ \diam B^n_w\leq C \sigma^n$.
\end{lem}

Let
$$
\OO\equiv \OO_F= \bigcap_{n=1}^\infty \bigcup_{w\in W^n} B^n_w.
$$
Let us also consider the {\it diadic group} $\displaystyle { W^\infty = \lim_{\longleftarrow} W^n}$.
The elements of $W^\infty$ are infinite sequences $(w_1 w_2\dots)$ of symbols $v\equiv 0$ and $c\equiv 1$
that can be also represented as formal power series
$$
   w \mapsto\sum_{k=0}^\infty w_{k+1} 2^k.
$$
The integers $\Z$ are embedded into $W^\infty$ as finite series.
The {\it adding machine} $p: W^\infty\ra W^\infty$ 
is the operation of adding $1$ in this group.   
The discussion above yields that 
 the map $F$ acts on the invariant Cantor set
$\OO$ as the dyadic adding machine 
(as in  the one-dimensional case, compare \cite{Mi}): 

\begin{cor}\label{adding machine}
The map $F|\OO$ is topologicaly conjugate to $p: W^\infty\ra W^\infty$. 
The conjugacy is given by the following homeomorphism $h: W^\infty\ra \OO$: 
$$
     h: w = (w_1w_2\dots) \mapsto \bigcap_{n=1}^\infty B^n_{w_1\dots w_n}. 
$$
 Furthermore,
$$
\text{HD}(\OO)\le \frac{\log 2}{\log \lambda}\le 0.73. 
$$
\end{cor}

We call $\OO$ the {\it critical Cantor set}%
\footnote{This Cantor set consists of the ``critical points'' of $F$.
More precisely, we will show in the forthcoming notes that generically $\OO$ is the set 
of  singularities of the unstable lamination of $F$.} 
of $F$. 
Let us finish this section with a remark on the dependence of this Cantor set on~$F$:

\begin{lem}\label{hol dependence}
    The critical Cantor set $\OO_F\subset \Om$ moves holomorphically as $F$ ranges
over $\II^c_\Om(\bar\eps)$.
\end{lem}

\begin{proof}
  Each  contraction $\Psi^n_w=\Psi^n_w(F)$, $w\in W^n$,   has a unique attracting fixed point $\alpha^n_w(F)$.
By the Implicit Function Theorem, this point depends holomorphically on $F$.

By Lemma \ref{adding machine}, 
any point of $\OO_F$ can be encoded as $\alpha^\infty_w(F)$, where
$w=(w_1,w_2\dots)\in W^\infty$. 
Lemma \ref{boxes shrink} implies that $\alpha^n_{w_1\dots w_n}(F) \to \alpha^\infty_w(F)$ as $n\to \infty$,
at an exponential rate uniform in $F$. 
Since uniform limits of holomorphic functions are holomorphic, $\alpha^\infty_w(F)$ depends holomorphically on $F$. 

Moreover, since the coding  $h: W^\infty\ra \OO_F$ is injective, 
$\alpha^\infty_w(F) \not= \alpha^\infty_v(F)$ if $v\not= w$, and the conclusion follows.
\end{proof}

\section{The average Jacobian}
\label{exponents}

In this section we consider the average Jacobian $b$ of an infinitely
renormalizable H\'enon-like map with respect to the unique invariant
measure supported on its critical Cantor set. It is shown that the
characteristic exponents of this measure are 0 and $\log b$ and that
$b$ is a natural parameter for infinitely renormalizable
maps.

\bigskip

We continue to consider infinitely renormalizable H\'enon-like maps
and assume, moreover, that they are diffeomorphisms. They are, however,
allowed to be complex. 
Lemma~\ref{contracting} and the standard distortion estimate imply:

\begin{lem}[Distortion Lemma]\label{distortion}
   There exist constants $C$ and  $\rho<1$ such that for any piece $B^n_w$
and for any $y,z\in B^n_w$, $w\in W^n$ the following holds:
$$
   \log \left| {\Jac F^k(y) \over \Jac F^k(z)} \right| \leq C\rho^n, \quad
k=1,2,\dots, 2^n.
$$
\end{lem}

Since $F|_\OO$ is the adding machine, it has a unique invariant
measure $\mu$. Let us consider the  average Jacobian with respect
to this measure:
  $$
      b =\exp \int \log \Jac F\, d\mu.
   $$

\begin{cor}\label{Jac}
For any piece $B^n_w$ and any point 
$z \in B^n_w$,
$$
     \Jac F^{2^n}(z) = b^{2^n}(1+O(\rho^n)),
$$
where $\rho$ is as in Lemma~\ref{distortion}.
\end{cor}

\begin{proof}
  Since
$$
    \int_{B^n_w}\log \Jac F^{2^n}\, d\mu =
            \int_\OO \log\Jac F\, d\mu = \log b,
$$
there exists a point $\zeta\in B^n_w$ such that
$$
       \log \Jac F^{2^n}(\zeta) = \log b/\mu(B^n_w)= 2^n \log b,
$$
and the assertion follows from Lemma~\ref{distortion}.
\end{proof}

The two characteristic exponents, $\chi_-\leq
\chi_0$, of the measure $\mu$ are given by  

\begin{thm}\label{char exp}
   The characteristic exponents of $\mu$ are $\chi_-= \log b$ and $\chi_0= 0$.
\end{thm}

\begin{proof}
Let $G_n$ be the $n$-th renormalization of $F$. This map is smoothly
conjugate to the restriction of $F^{2^n}$ to the piece
$B^n_{v^n}$. Let $\mu_n$ be the normalized restriction of $\mu$ to
$B^n_{v^n}$, and let $\nu_n$ be the invariant measure on the critical Cantor set of
$G_n$. Note that these two measures are preserved by the
conjugacy. Then 
$$ 2^n \chi_0 = \chi_0 (F^{2^n}| B^n_{v^n},\, \mu_n) =
\chi_0(G_n , \nu_n) \leq
\int \log \|DG_n\|\,  d\nu_n \leq C,
$$
since the maps $G_n$ have uniformly bounded $C^1$-norms.

Hence $\chi_0\leq 0$. If $\chi_0<0$, both characteristic exponents of
$F$ would be negative and it would then follow from the Pesin theory that
$\mu$ is supported on a periodic cycle%
\footnote{Indeed, in this case the Pesin local stable manifold $W=W^s_{\loc}(x)$ (see e.g., \cite{PS})
  of a typical point $x\in \OO$ would be a neighborhood of $x$. Then for some big $n$,
   $f^n$ would be a contracting map of $W$ into itself, and the $\orb x$ would converge to an
    attracting cycle.}    
which is clearly not the case. Hence
$\chi_0=0$. The formula for the other exponent now follows from the
relation $\chi_0+ \chi_- = \log b.$
\end{proof}

\msk Let us now take a look at the dependence of the average
Jacobian on parameters. Consider a holomorphic one-parameter
family of complex H\'enon-like maps $F_t\in \II^c_\Om(\bar\eps)$,
$$
    F_t \colon (x,y)\mapsto (f(x)-t\, \eps_t(x,y),\, x),\quad |t|<r,\ (x,y)\in \Om,
$$
 such that
\begin{enumerate}[(i)]
\item $F_t$ are real for real $t$;

\item $\eps_t(x,y)= \gamma(x,y)\psi_t(x,y)$, where
$\psi_t(x,y) = 1+ O(t)$;

\item $\di \gamma / \di y >0$ on $B$ and $\di \gamma /
\di y\not=0 $ on $\Om$.
\end{enumerate}

Let us consider the {\it complex Jacobian},
$$
  \Jac^c F_t=\det DF_t = t{\di \eps_t\over \di y} = t {\di\gamma\over \di y}
+O(t^2).
$$
By property (iii), it  does not vanish for sufficiently small $r$,
and hence $F_t$ are complex diffeomorphisms. Moreover, for real
$t$, they preserve orientation of $B$.

\comm{ Let us consider a parameter sector
$$
  S= \{t\in \C: 0< |t| <r, \ |\arg t|< \theta\ {\mathrm{for some}}\ \theta<\pi
\}.
$$
The function  $(t,x,y)\mapsto \log \Jac^c F_t(x,y)$ admits a
holomorphic extension to $S \times \Om$ which is positive for
positive $t$. Hence for $t\in S$ we can define the {\it complex
average Jacobian}
$$
   b^c(F_t) = \exp \int \log \Jac F_t d\mu_t
$$
where $\mu_t$ is the $F_t$-invariant measure on the critical Cantor set
$\OO_{F_t}$. Since by Lemma~\ref{hol dependence}, the measure
$\mu_t$ depends holomorphically on $t$, we conclude: }

\begin{lem}\label{b-dependence}
For sufficiently small $r>0$, the average Jacobian $b_t\equiv b(F_t)$,
$t \in (0,r)$, admits a holomorphic extension to the complex disk
$\D_r$. Moreover,
\[
    b'(0) = \exp \int_{O(f)} \log {\di \gamma \over \di y}\, d\mu\neq 0.
\]
\end{lem}

\begin{proof}
We can define the average complex
Jacobian by the following explicit formula: 
$$
       b^c(F_t) = \exp \int_{\OO_t} \log \Jac^c F_t\, d\mu_t= 
$$
$$ = t \exp \int_{\OO_t} \log {\di\gamma \over \di y}\, d\mu_t \, \cdot\,
               \exp \int_{\OO_t} \log \psi_t(x,y)\, d\mu_t ,
$$
where $\mu_t$ is the $F_t$-invariant measure on the critical Cantor set
$\OO_t= \OO_{F_t}$. Since $\psi_t = 1+ O(t)$, there is a well
defined holomorphic branch of  $\log \psi_t(x,y)$ on the domain
$\D_r \times \Om$ which  is positive on $(-r,r)\times B$. Since
by Lemma~\ref{hol dependence} the Cantor set $\OO_t$ moves
holomorphically with $t$, the two integrals on the right-hand side
of the formula above depend holomorphically on $t$. Since the
second factor in that product goes to 1 as $t\to 0$, the result follows.
\end{proof}

Thus, in the H\'enon-like families as above, the average Jacobian
$b$ can be used (consistently with the common intuition) as a
holomorphic parameter that measures the distance to the reference
unimodal map.

\section{Universality around the tip}\label{univ}

This section is central in our paper.
We  prove here that the renormalizations of H\'enon-like maps
near the tip have the following shape: 
%
 \[ R^nF=(f_n - b^{2^n}a(x)\, y\, (1+O(\rho^n)),\,x),  \]
where $a(x)$ is a {\it universal} function associated with the unimodal fixed point $f_*$. 
To establish this Universality Law, we study closely the
Renormalization Microscope constructed in Section~\ref{Cantor}.
Lemma~\ref{APPsi}, Lemma~\ref{ustar}, and Corollary~\ref{tilt}
are the main technical results of this section; they 
quantify the {\em tilting} phenomenon mentioned earlier. These
lemmas will also be crucial in the next sections when the non-rigidity and the
existence of critical Cantor sets with unbounded geometry is established. 

\subsection{Some universal one-dimensional  functions}\label{1D universal f-s} 

Recall that\break  $f_*\colon I\to I$ stands for the one-dimensional
renormalization fixed point 
normalized so that $f_*(c_*)=1$ and $f_*^2(c_*)=-1$, where 
$c_*\in I$ is the  critical point of $f_*$. 
We let $J_c^* =[-1, f_*^4(c_*)]$ be the smallest renormalization interval of $f_*$, 
and   we let $s\colon J_c^*\ra  I$ be the orientation reversing affine rescaling.
The smallest renormalization interval around the critical value is denoted by
 $J_v^*= f_*(J^*_c)=[f_*^3(c_*),1]$.
Then $s\circ f_*: J_v^*\ra [-1,1] $ is an expanding  diffeomorphism.
Let us consider the inverse contraction  
\[
    g_*\colon I\to J^*_v, \quad g_*=f_*^{-1}\circ s^{-1},
\]
where $f_*^{-1}$ stands for the branch of the inverse map that maps $J_c^*$ onto  $J_v^*$. 
The function $g_*$ is the non-affine branch of the so called ``presentation function''
(see~\cite{BMT} and references therein).
Note that $1$ is the unique fixed point of $g_*$.

 Let $J_c^*(n) \subset I$ be the smallest periodic interval of period $2^n$ that contains $c_*$ and $J_v^*(n) \subset I$ be the smallest periodic interval of period $2^n$ that contains $1$.

Let
$G_*^n\colon I\to I$ be the diffeomorphism obtained by rescaling
affinely the image of $g_*^n$. The fact that $g_*$ is a contraction implies
that the following limit exists 
\[
u_*=\lim_{n\to \infty} G_*^n \colon I\to I,  
\]
where the convergence is exponential in the $C^3$-topology.
In fact, this function linearizes $g_*$ near the attracting fixed point $1$ 
(see, e.g., \cite[Theorem 8.2]{M}).

\begin{lem} \label{ustarfstar} For every $n\ge 1$
\begin{enumerate}
\item [\rm (1)] $J^*_v(n)=g_*^n(I)$,

\item [\rm (2)] $R^n_vf_*= G^n_*\circ f_*\circ (G^n_*)^{-1}$.
\end{enumerate}
Moreover,
\begin{enumerate}
\item [\rm (3)]
$
u_*\circ f^*=f_*\circ u_*.
$
\end{enumerate}
\end{lem}

\begin{proof} The proof of the first two items is by induction. 
Notice that the definition of $g_*$ implies directly
$$
f_*^{2}|J^*_v= g_* \circ f_* \circ (g_*)^{-1}.
$$
Let $h_n:I\to J^*_v(n)$ be the conjugation between the two infinitely 
renormalizable maps $f_*^{2^n}|J^*_v(n)$ and 
$f_*$,
$$
f_*^{2^n}|J^*_v(n)= h_n \circ f_* \circ (h_n)^{-1}.
$$
Note, $h_1=g_*$. A calculation shows,
$$
h_{n+1}=h_n\circ g_*.
$$
To do this calculation,  first notice that
$$
J^*_v(n+1)=h_n(J^*_v).
$$
Hence,
$$
\begin{aligned}
f_*^{2^{n+1}}|J^*_v(n+1)&=f_*^{2^{n}}|J^*_v(n) \circ f_*^{2^{n}}|J^*_v(n+1)\\
                      &=h_n \circ f_*^2|J^*_v  \circ (h_n)^{-1}\\
                      &=(h_n \circ g_*) \circ f_*^2 
                        \circ (h_n \circ g_*)^{-1}.
\end{aligned}
$$
Now, $R^n_vf_*$ is obtained by rescaling $f_*^{2^n}|J^*_v(n)$. In particular,
$$
R^n_vf_*= G^n_*\circ f_*\circ (G^n_*)^{-1}.
$$
This finishes the proof of item $(1)$ and $(2)$. 
The convergence of the sequence $G^n_*$ to $u_*$ implies that 
$R^n_vf_*$ converges. The limit has to be the unique fixed point $f^*$ 
of $R_v$. This finishes the proof of $(3)$.
\end{proof}

Notice that $|J^*_c(n)|=\sigma^n$ and 
$
f_*(J^*_c(n))=J^*_v(n)=g^n_*(I)$.
Hence,

\begin{cor}\label{Dgstar}
$\displaystyle
\frac{dg_*}{dx}(1)=\sigma^2.
$
\end{cor}

Along with $u_*$, we consider its rescaling 
$$
v_*: I\ra \R, \quad v_*(x)=\frac{1}{u'_*(1)}(u_*(x)-1)+1,
$$
normalized so that $v_*(1)=1$ and $\displaystyle \frac{d v_*}{dx} (1)=1$.

\begin{lem}\label{convergence to g-star}
  Let $\rho\in (0,1)$, $C>0$.  
Let us consider a sequence of smooth functions $g_k: I\ra I$, $k=1,\dots, n$,
 such that 
$\| g_k - g_*\|_{C^3}\leq C \rho^k$. Let
$g^n_k=g_k\circ\dots\circ g_n$, 
and let         
$G^n_k= a^n_k\circ g^n_k: I\ra I$, where $a^n_k$ is the affine rescaling of 
$\Im g^n_k$
 to $I$. 
Then  $\|G^n_k - G_*^k\|_{C^2}\leq C_1 \rho^{n-k}$, where $C_1$ depends only 
on $\rho$ and $C$.
\end{lem}   
 
\begin{proof}   
Let $I_0=I$ and $I_j=[x_j,y_j]\subset I$ such that $g_j(I_j)=I_{j-1}$. 
Rescale affinely the domain and image of $g_j\colon I_j\to I_{j-1}$
and denote the normalized diffeomorphism by $h_j\colon [-1,1]\to [-1,1]$.
Let
$$
I_j^*=[x_j^*,1]=g_*^{n-j}([-1,1])
$$
and rescale the domain and image of $g_*\colon  I_j^*\to I_{j-1}^*$ and
denote the normalized diffeomorphism by $h_j^*\colon [-1,1]\to
[-1,1]$. Note that
$$
h_{k}^*\circ h_{k+1}^*\circ \cdots \circ h_n^* \to
u_*,
$$
where the convergence in the $C^2$ topology is exponential in $n-k$. 
In the following estimates we 
will use a uniform constant $\rho<1$
for exponential estimates.
Let $\Delta x_j=x_j-x^*_j$ and $\Delta y_j=1- y_j$  . Then
$$
x_{j-1}=g_*(x^*_j)+g_*'(z) \cdot \Delta x_j +O(\rho^j).
$$
Hence, using a similar argument for  $\Delta y_j$,
$$
|\Delta x_j|, |\Delta y_j|= O(\rho^j).
$$
Because, $g_j$ and $g_*$ are contractions we have
$$
|I_j|, |I^*_j|=O(\rho^{n-j}).
$$
We will represent a diffeomorphism  $\phi:I\to J$ by its nonlinearity
$$
\eta_\phi=\frac{D^2\phi}{D\phi}.
$$
Let $\eta_j$ and $\eta^*$ be the nonlinearities of $g_j$ and $g_*$.
Notice that
$$
\|\eta_j-\eta^*\|_{C^1}=O(\rho^j).
$$
Furthermore, let  $\mathbb{I}_j:[-1,1]\to I_j$ and 
$\mathbb{I}^*_j:[-1,1]\to I^*_j$
be the affine orientation preserving rescalings. Using this notation 
$$
\eta_j(\mathbb{I}_j(x))=\eta^*(\mathbb{I}^*_j(x))+
D\eta^*(z)\cdot \left(\mathbb{I}_j(x)-\mathbb{I}^*_j(x)\right)+     O(\rho^j),
$$
for some $z\in [\mathbb{I}_j(x),\mathbb{I}^*_j(x)]$.
Hence,
$$
\eta_j(\mathbb{I}_j(x))=\eta^*(\mathbb{I}^*_j(x))+
     O(\rho^j).
$$
The nonlinearities of $h_j$ and $h^*_j$ are given by
$$
\eta_{h_j}=|I_j|\cdot \eta_j(\mathbb{I}_j),
$$
and similarly
$$
\eta_{h^*_j}=|I^*_j|\cdot \eta^*(\mathbb{I}^*_j).
$$
Now
$$
|\eta_{h_j}(x)-\eta_{h^*_j}(x)|=O((|I_j|-|I^*_j|)+ \rho^j\cdot |I^*_j|).
$$
Hence
$$
|\eta_{h_j}(x)-\eta_{h^*_j}(x)| =
\left\{
\begin{array}{ccc }
O(\rho^{n-j}) & : & j\leq (n+k)/2 
\\ 
O(\rho^{j}) & : & j> (n+k)/2.
\end{array} \right.
$$
It follows that
$$
\sum_{j=k}^n \|\eta_{h_j}- \eta_{h^*_j}\|_{C^0}=O(\rho^{n-k}).
$$
Note that we can estimate $\|\eta_{h_j}\|_{C^1}$ by using
$$
D\eta_{h_j}=|I_j|^2 D\eta_{h_j}(\mathbb{I}_j).
$$
The resulting estimate allows to use
a reshuffling argument,   
see Appendix , Lemma~\ref{shufflem}, which finishes the
proof of the Lemma.
\end{proof}

\subsection{Asymptotics of  the $\Psi$-functions}\label{Psi-functions}
Fix an infinitely renormalizable H\'enon-like map $F\in \II_\Om(\bar\eps)$ 
to which we can apply Theorem~\ref{convergence}.
For such an $F$, we have a well defined {\it tip}: 
$$
\tau\equiv \tau(F)=\bigcap_{n\ge 0} B^n_{v^n},
$$
where the pieces $B^n_w$ are introduced in \S \ref{pieces}. 
Let us consider the tips of the renormalizations, $\tau_k=\tau(R^k F)$.
To simplify the notations, we will translate these tips to the origin
by letting
$$
     \Psi_k\equiv \Psi_k^{k+1} =  \Psi^1_v (R^k F)\,  (z + \tau_{k+1}) - \tau_k.  
$$ 

Denote the derivative of the maps $\Psi_k $ at $0$ by $D_k\equiv D_k^{k+1}$
and decompose it into the unipotent and diagonal factors: 
 \begin{equation}\label{Dk}
 D_k= \left(
\begin{array}{cc}
1 & t_k\\
0 & 1
\end{array}\right)
\left(
\begin{array}{cc}
\alpha_k & 0\\
0 & \beta_k
\end{array}\right). 
\end{equation}
Let us factor this derivative out from $\Psi_k$: 
$$
\Psi_k = D_k \circ (\id + {\bf s}_k),
$$
where ${\bf s}_k(z) = (s_k(z) , 0) =  O(|z|^2)$ near 0.
The convergence Theorem~\ref{convergence}
and the explicit expression for the $\Psi$-functions
(see (\ref{H}) and \S \ref{branches}) imply:

\begin{lem} \label{smalls}
There exists $\rho<1$ such that for $k\in \Z_+$ the following estimates hold:


\item [\rm (1)] $\displaystyle 
  \alpha_k=\sigma^2 \cdot (1+O(\rho^k)),\quad \beta_k=-\sigma \cdot (1+O(\rho^k)), \quad
         t_k=O(\bar\eps^{2^k}); 
$

\item [\rm (2)]
$\displaystyle | \di_x s_k | = O(1),\quad  |\di_y s_k| = O(\bar\eps^{2^k}); $  

\item[\rm (3)]
$\displaystyle
  |\partial^2_{xx} s_k |= O(1),\quad
  |\partial^2_{xy} s_k |= O(\bar{\eps}^{2^k}) ,\quad 
  |\partial^2_{yy} s_k| =O(\bar{\eps}^{2^k}).
$

\end{lem}

Note that since all the maps under consideration are holomorphic, 
the bounds on their derivatives follow from  the bounds on the maps themselves.

Let now
$$
   \Psi_k^n = \Psi_k\circ\dots \circ \Psi_{n-1}, \quad B_k^n= \Im \Psi_k^n.
$$
Since by Lemma \ref{contracting}
$$
{\diam}(B_k^n)=O(\sigma^{n-k})\quad {\rm for} \quad  \quad k<n,
$$
we conclude:

\begin{cor}\label{second derivatives} 
 Let $k<n$. For $z \in B_{k+1}^n$ we have: 
$$
\left| \partial_x s_k(z)\right|   =O (\sigma^{n-k}), \quad 
\left| \partial_y s_k(z)  \right| =
      O(\bar{\eps}^{2^k}\cdot\sigma^{n-k}).
$$
\end{cor}

Let us now consider the derivatives of the maps $\Psi^n_k$ at the origin:  
$$
D_k^n=D_k\circ D_{k+1}\circ \cdots D_{n-1}.
$$
Since the unipotent matrices form a normal subgroup in the group of
upper-triangular matrices, we can reshuffle this composition and obtain:

 \begin{equation}\label{reshuffling}
 D_k^n= 
 \left(
\begin{array}{cc}
1 & t_k\\
0 & 1
\end{array}\right)
\left(
\begin{array}{cc}
(\si^2)^{n-k} & 0\\
0 & (-\si)^{n-k}
\end{array}\right) (1+O(\rho^k)). 
\end{equation}  
Factoring the derivatives  $D_k^n$ out from $\Psi_k^n$, we obtain: 
\begin{equation}\label{factoring}
\Psi_k^n = D_k^n \circ (\id + {\bf S}_k^n),
\end{equation}
where ${\bf S}_k^n (z) = (S_k^n(z), 0) = O(|z|^2)$ near 0. 

\begin{lem}\label{APPsi} For $k<n$,  we have:
\begin{enumerate}

\item [{\rm (1)}]
$ \displaystyle
   | \di_x S^n_k | = O(1),\quad  |\di_y S^n_k | = O(\bar{\eps}^{2^k}); 
$

\item [{\rm (2)}]
$ \displaystyle
|\partial^2_{xx} S^n_k| =O(1), \quad
|\partial^2_{yy} S^n_k| =O(\bar{\eps}^{2^k}), \quad
|\partial^2_{xy} S^n_k| =O(\bar{\eps}^{2^k}\, \si^{n-k}).
$
\end{enumerate}

\end{lem}

\begin{proof} 
Let
$$
z^n_{k+1}=
\left(
\begin{array}{c}
x^n_{k+1}\\
y^n_{k+1}
\end{array}\right)
=\Psi_{k+1}^n (z)
$$ 
By (\ref{reshuffling}) and (\ref{factoring}),  
\begin{eqnarray*}
 x^n_{k+1} &=& K_1\, (\sigma^2)^{n-k-1}\,
(x+S^n_{k+1}(x,y))+ K_2\, t_k
\,  (-\sigma)^{n-k-1}\, y \\
y_n^{k+1} &=& K_3\,  (-\sigma)^{n-k-1}\,  y ,
\end{eqnarray*}
where $K_i= K_i(k,n) = O(1) $
(and the constants $K_i$ below have the same meaning).

Moreover, since 
$$
D^n_k\circ (\id +{\bf S}^n_k) =  \Psi_k^n = \Psi_k\circ \Psi^n_{k+1}= 
  D_k \circ   (\id + {\bf s}_k)  \circ \Psi_{k+1}^n  =
$$
$$
 D_k^n \circ(\id + {\bf S}^n_{k+1})  +    D_k \circ  {\bf s}_k \circ \Psi^n_{k+1},
$$
we obtain:
$$
S^n_k(z)= S^n_{k+1}(z)+ K_4\, (\la^2)^{n-k-1} \, s_k (z^n_{k+1})
$$
(recall that $\la= \si^{-1}$). 

\medskip

The proof proceeds by relating the partial derivatives of $S^n_k$ to the
derivatives of $ S^n_{k+1}$ and $s_k$. For instance, by differentiating the last equation
taking into account the above expressions for $x^n_{k+1}$ and $y^n_{k+1}$, we obtain: 
$$
  \frac{\partial S^n_k}{\partial y} = 
  (1+K_5\,  \frac{\partial s_k}{\partial x} ) \, \frac{\partial S_n^{k+1}}
{\partial y}  + 
  K_6 \,  t_k\,  (-\la )^{n-k-1}\, \frac{\partial s_k}{\partial x}  +   
  K_7\, (-\la)^{n-k-1} \, \frac{\partial s_k}{\partial y} ,
$$
where the partial derivatives of $s_k$ are computed at
$z^n_{k+1}$. 
Now Corollary~\ref{second derivatives} implies 
$$
\left|\frac{\partial S^n_k}{\partial y}\right|\le (1+O(\rho^{n-k}))\, 
\left|\frac{\partial S^n_{k+1}}{\partial y}\right|+
 C \, \bar{\eps}^{2^k},
$$
and hence for all $k<n$, 
$$
\left|\frac{\partial S^n_k}{\partial y}\right|\le  C \,
\bar{\eps}^{2^k},
$$
as was asserted.
The bound for  $\di S^n_k/ \di y$ is obtained in a similar way. 

Since the functions $S^n_k$ are holomorphic and defined on a fixed domain,
the first two  bounds on the second derivatives follow. 
However, the bound on the mixed  derivative does not follow from this general reasoning. 
Differentiating $\di S^n_k/ \di y$ (taking into account the expressions for $x^n_{k+1}$ and $y^n_{k+1}$), 
we obtain:

$$
\begin{aligned}
\frac{\partial^2 S^n_k}{\partial xy}= &
\left(1+K_5 \frac{\partial s_k}{\partial x}\right)\, \frac{\partial^2 S^n_{k+1}}{\partial xy}+\\
& \left(1+\frac{\partial S_{k+1}^n}{\partial x}\right)\, (\sigma^2)^{n-k-1} \, \frac{\partial^2 s_k}{\partial x^2}
         \left(K_8\, \frac{\partial S_{k+1}^n}{\partial y}+  K_9\, t_k \la^{n-k-1}  \right)+\\
& K_{10}\left(1+\frac{\partial S_{k+1}^n}{\partial x}\right)\, 
  (-\sigma)^{n-k-1}\,  \frac{\partial^2 s_k}{\partial xy}, 
\end{aligned}
$$
where the partial derivatives of $s_k$ are calculated at
$x^n_{k+1}$.
Using Corollary~\ref{second derivatives} and the previous estimates on the first
partial derivatives of $S_k^n$, we obtain
$$
\left|\frac{\partial^2 S^n_k}{\partial xy}\right|\le
(1+O(\rho^{n-k}))\cdot \left|\frac{\partial^2 S^n_{k+1}}{\partial xy}\right|+
     C \cdot \bar{\eps}^{2^k}\cdot \sigma^{n-k} .
$$
Hence, 
$$
\left|\frac{\partial^2 S^n_k}{\partial xy}\right|\le  C \cdot
\bar{\eps}^{2^k}\cdot \sigma^{n-k} .
$$
\end{proof}




\comm{**************************
Theorem \ref{convergence} implies the exponential
convergence
\[
\Psi^{k,k+1} \left(
\begin{array}{c}
x\\
y
\end{array}\right)
\to \left(
\begin{array}{c}
g_*(x+1)-1\\
-\sigma\cdot y
\end{array}\right)
\quad\text{ as } k\to\infty.
\]
*********************************}

We are now ready to describe the asymptotical behavior of the $\Psi$-functions
using the universal one-dimensional functions from  \S \ref{1D universal f-s}. 
Let us normalize the function $v_*$ so that it fixes $0$ rather than $1$: 
$$
   {\bf v}_*(x) = v_*(x+1)-1.
$$

\begin{lem}\label{ustar}
 There exists $\rho<1$ such that  for all $k<n$ and $y\in I$,
$$
\left|\id + S^n_k (\cdot,y)-{\bf v}_*(\cdot)\right|=O(\bar{\epsilon}^{2^k}\cdot y+
\rho^{n-k})
$$
and
$$
  \left| 1+ \frac{\partial S_k^n}{\partial x}(\cdot,y)-
\frac{\partial {\bf v}_*}{\partial x}(\cdot)\right|=O(\rho^{n-k}).
$$
\end{lem}

\begin{proof} By  Lemma \ref{APPsi}, 
$$
\left|\frac{\partial^2 S_k^n}{\partial yx} \right|=O(\bar\eps^{2^k}\sigma^{n-k})=O(\sigma^{n-k})
$$
and
$$
\left|\frac{\partial S_k^n}{\partial y} \right|=O(\bar{\epsilon}^{2^k}).
$$
Hence it is enough to verify the desired convergence on the horizontal section passing through the tip: 
$$
\dist_{C^1}(\id + S_k^n (\cdot, 0 ),\, {\bf v}_*(\cdot))=O(\rho^{n-k}).
$$

Let us normalize $g_*$ so that $0$ becomes its fixed point with $1$ as 
multiplier:   
$$  
          \Bg_*(x) = \frac{g_*(x+1)-1}{g_*'(1)}.
$$   
Now, $\id+S^n_k(\cdot, 0)$ is the rescaling of $\Psi^n_k(\cdot, 0)$ 
normalized so 
that
the fixed point $0$ has  multiplier 1. By Theorem \ref{convergence}, 
$$
    \dist_{C^3}(\id+ s_k(\cdot, 0), \Bg_*(\cdot)) = O(\rho^k).
$$
Hence, by  Lemma \ref{convergence to g-star},  
$$
    \dist_{C^1}(\id + S^n_k(\cdot,0), \Bg_*^{n-k}(\cdot)) = O(\rho^{n-k}) .
$$ 
Since $\Bg^n \to {\bf v}_*$ exponentially fast, the conclusion follows. 
\end{proof}

\begin{prop}\label{limit}
 There exists a coefficient $a_F\in \mathbb{R}$ and an absolute constant $\rho\in (0,1)$ such that 
$$
\left|(x + S^n_0(x,y))-({\bf v}_*(x)+a_F y^2)\right|=O(\rho^n).
$$
\end{prop}

\begin{proof} The image  of the vertical interval  $y\mapsto (0,y)$
under the map $\id + \BS_0^n$ is the graph of a function $w_n:I\to \mathbb{R}$ defined by
$$
w_n(y)=S_0^n(0,y).
$$
By the second part of  Lemma~\ref{ustar} we have: 
$$
\left|(x + S^n_0(x,y))-({\bf v}_*(x) + w_n(y) )\right|
=O(\rho^n).
$$
Let us show that the functions $w_n$ converge to a parabola. 
The identity
$$
   D_0^{n+1}\circ (\id + \BS^{n+1}_0) = \Psi^{n+1}_0=\Psi^n_0\circ \Psi_n= D^n_0\circ(\id +\BS^n_0)\circ D_n\circ (\id+{\bf s}_n),
$$
implies
$$
   \BS^{n+1}_0= {\bf s}_n+ D_n^{-1}\circ \BS^n_0\circ D_n\circ (\id+{\bf S}_n),
$$
so that
\begin{equation}\label{w sub n}
w_{n+1}(y)=s_n(0,y)+ \frac{1}{\alpha_n} S_0^n(\alpha_n s_n(0,y)+\beta_n t_n y,\, \beta_n y),
\end{equation}
where $\alpha_n, \beta_n$ and $t_k$ are the entries of $D_n$, see 
equation~(\ref{Dk}).
The estimate of $\di_y s_n$ from  Lemma \ref{smalls} implies:
\begin{equation}\label{s sub n}
s_n(0,y)=e_n y^2 + O(\bar\eps^{2^n} y^3),
\end{equation}
where $e_n=O(\bar{\epsilon}^{2^n})$.
The estimate of $\di_{xy}^2 S^n_0$ from Lemma \ref {APPsi} implies:
$$
  \frac{\di S_0^n}{\di x}(0,y) = O(\bar\eps^{2^n} y).
$$
Hence
$$
S_0^n(\alpha_n s_n(0,y)+\beta_n t_n y,\, \beta_n y)
$$
$$
= S^n_0(0,\beta_n y)+ 
      \frac{\di S^n_0}{\di x} (0, \beta_n y) (\alpha_n s_n(0,y)+\beta_n t_n y)+O(\bar\eps^{2^n} y^3)
$$
$$
   = S^n_0(0, \beta_n y) + q_n y^2 + O(\bar\eps^{2^n} y^3 )= w_n(\beta_n y) + q_n y^2 + O(\bar\eps^{2^n} y^3 ),
$$
where $q_n= O(\bar\eps^{2^n})$.
Incorporating this and (\ref{s sub n}) into (\ref{w sub n}), we obtain: 
$$
w_{n+1}(y)= \frac{1}{\alpha_n} w_n(\beta_n y)+
             c_n y^2 + O(\bar\eps^{2^n} y^3),
$$
where $c_n=O(\bar{\epsilon}^{2^n})$. 
Writing $w_n$ in the form
$$
w_n(y)=a_ny^2 +A_n(y) y^3, 
$$
we obtain: 
$$
a_{n+1}=\frac{\beta_n^2}{\alpha_n} a_n +c_n
$$
and
$$
\|A_{n+1}\|\le \frac{|\beta_n|^3}{\alpha_n}\|A_n\|+ O(\bar\eps^{2^n}).
$$
Now  the first item of Lemma~\ref{smalls} 
implies that $a_n\to a_F$ and $\|A_n\|\to 0$ exponentially fast.
\end{proof}

\comm{
\
\begin{proof}

$$
\dist_{C^1}(h_k\circ h_{k+1}\circ \cdots \circ h_n,
h^*_k\circ h^*_{k+1}\circ \cdots \circ h^*_n)=O(\rho^{n-k}).
$$
Notice that $u^k_n$ is an affine rescaling of
 $h_k\circ h_{k+1}\circ \cdots \circ h_n$. We can easily recover this 
rescaling by observing that 
$$
\frac{\partial u^k_n}{\partial x}(\tau_{n,1},\tau_{n,2})=1,
$$
see Lemma \ref{APPsi}. In particular
$$
u^k_n(x)=\frac{h_k\circ h_{k+1}\circ \cdots \circ h_n(x)-h_k\circ h_{k+1}\circ \cdots \circ h_n(\tau_{n,1})}
{(h_k\circ h_{k+1}\circ \cdots \circ h_n)'(\tau_{n,1})}+\tau_{k,1}.
$$
This then implies
$$
\dist(u^k_n(\cdot,\tau_{n,2}), v_*)=O(\rho^{n-k})
$$
for some $\rho<1$.
\end{proof}
}

\subsection{Universality}

\comm{********************************
 Note first that since the boxes $B^n_v$ 
shrink uniformly exponentially in $n$, the standard distortion estimates yield
(compare Lemma \ref{distortion}):

\begin{lem}\label{Jac Psi}
  The maps $\Psi^n $  have uniformly bounded  Jacobian distortion. 
\end{lem}

\begin{lem}\label{de}
Let $F_n := R^n F =  (f_n -\eps_n, x)$.  Then 
$$ {\di\eps_n \over \di y}
\asymp b^{2^n}. $$ 
\end{lem}

 According to the chain rule,

\begin{equation}\label{chain rule}
\begin{aligned}
{\di\eps_n \over \di y}(z)& = \Jac F_n(z) \\
&= \Jac F^{2^n}(\Psi^n(z))\cdot\frac{\Jac \Psi^n (z)}{\Jac\Psi^n(F_n (z))}.
\end{aligned}
\end{equation}
Since $\Psi_n$ has a bounded Jacobian distortion by Lemma~\ref{Jac Psi},
the conclusion follows from Corollary~\ref{Jac}.  
\end{proof}

\begin{cor} The numbers $t_k$ defined by equation~(\ref{tk}) satisfy
$$
t_k\asymp  b^{2^k}.
$$
\qed
\end{cor}
****************************************************}

We are ready to  prove the main positive result of this paper: 

\begin{thm}[Universality]\label{universality}
For any $F\in \II_\Om(\bar\eps)$ with sufficiently small $\bar\eps$, we have: 
\[
    R^n F = (f_n(x) -\,  b^{2^n}\, a(x)\, y\, (1+ O(\rho^n)), \ x\, ),
\]
where $f_n\to f_*$ exponentially fast, $b$ is the average Jacobian, $\rho\in (0,1)$,
and $a(x)$ is a universal  function. Moreover, $a$ is analytic and 
positive.   
\end{thm}

\begin{proof} 
Let $F_n\equiv R^n F= (f_n-\eps_n,x)$. 
The function  $\Psi^n\equiv \Psi^n_v$  conjugates
the renormalization $F_n$ to the iterate $F^{2^n}$ on the piece $B^n\equiv B^n_v$. 
(Here $\Psi^n$ is the  original $\Psi$-function rather than the normalized one, $\Psi^n_0$.)
According to the chain rule,

\begin{equation}\label{chain rule}
\begin{aligned}
 \di_y \eps_n (z)& = \Jac F_n(z) 
  = \Jac F^{2^n}(\Psi^n(z))\, \frac{\Jac \Psi^n (z)}{\Jac\Psi^n(F_n z)} \\
&=   b^{2^n}  \,  \frac{\Jac \Psi^n (z)}{\Jac\Psi^n(F_n z)} \, (1+ O(\rho^n)) ,
\end{aligned}
\end{equation}
where the last equality follows from Lemma~\ref{Jac}.

Let $D^n\equiv D^n_0$, $\BS^n\equiv \BS^n_0$,  $S^n\equiv S^n_0$.
Let us consider affine maps $T^n: z\mapsto z-\tau_n$ and $L^n: z\mapsto (D^n)^{-1} (z-\tau)$
as local charts on $B$ and $B^n$ respectively. Various maps presented in these local charts 
will be written in the boldface, so that
$$
  \BF_n = T^n\circ F_n\circ (T^n)^{-1}, \quad \BPsi^n\equiv \id + \BS^n =  L^n\circ \Psi^n\circ(T^n)^{-1}. 
$$

Since affine maps do not distort the Jacobian, we have:  
\begin{equation}\label{JacPsi}
  \frac{ \Jac\Psi^n(z)}{\Jac \Psi^n(F_n z)} = \frac{\Jac \BPsi^n(\Bz)}{\Jac \BPsi^n (\BF_n\Bz)} = \frac{1+ \di_x S^n(\Bz)}{1+\di_x S^n(\BF_n \Bz)}, 
\end{equation}
where $\Bz= Tz$. 

By Lemma \ref{ustar}, 
\begin{equation}\label{Bv-star}
1+ \di_x S^n \to \Bv_*'
\end{equation}
 exponentially fast. 
By Theorem \ref{convergence}, $\tau_n\to \tau_\infty \equiv (c_*,1)$ exponentially fast, so that
$T_n$  converges exponentially  to the translation $T^\infty: z\mapsto z-\tau_\infty$.
Applying Theorem \ref{convergence}  once again, we conclude that 
$\BF_n \to (\Bf_*, x)$ exponentially fast, where $\Bf_*(x) = f_*(x+1)-1$.
Putting this together with (\ref{JacPsi}) and (\ref{Bv-star}), we conclude:
$$
    \frac{ \Jac\Psi^n(z)}{\Jac \Psi^n(F_n z)} \to \frac{\Bv_*'(\Bx)}{\Bv_*'(\Bf_*(\Bx))} =  \frac{v_*'(x)}{v_*'(f_*(x))} \equiv a(x), 
$$
where $z=(x,y)$, $\Bx=x-1$, and convergence is exponential. Since $v_*$ is an analytic diffeomorphism,
the function $a(x)$ is analytic and non-vanishing. 

Plugging the last formula into (\ref{chain rule}), we obtain:
$$
   \di_y \eps_n(z)  =  b^{2^n}\, a(x) \, (1+ O(\rho^n)) .
$$
Integration of this formula yields: 
$$
    \eps_n(x,y) = c_n(x) + b^{2^n}\, a(x) \, y\, (1+ O(\rho^n)),
$$
and since $\|c_n(x)\| $ is super-exponentially small,
 it can be incorporated into the unimodal term $f_n(x)$.
\end{proof}

\begin{cor}\label{tilt} The numbers $t_k$ defined by 
equation~(\ref{reshuffling}) satisfy$$
t_k\asymp  -b^{2^k}.
$$
\end{cor}

\begin{proof} Consider $\Psi_k=(\Lambda_k\circ H_k)^{-1}$, where
$\Lambda_k$ and $H_k$ are used to define $R^{k+1}F$. Recall
$$
\Lambda_k(x,y)=\left(
\begin{array}{c}
s_k(x) \\
s_k(y)
\end{array}\right)
$$
and
$$
H_k(x,y)=\left(
\begin{array}{c}
f_k(x)-\epsilon_k(x,y) \\
y
\end{array}\right)
$$
where $s_k$ is an orientation reversing affine map with $s\asymp -1$ as 
derivative. Then
$$
D_k^{-1}=D\Lambda_k\circ DH_k =
\left(
\begin{array}{cc}
\cdot  & -s\di_y \eps_k(\tau_k)\\
0 & \cdot
\end{array}
\right).
$$
The representation of $D_k$ from (\ref{Dk}) gives
$$
 \left(
\begin{array}{cc}
1 & -t_k\\
0 & 1
\end{array}\right)
=
\left(
\begin{array}{cc}
\alpha_k & 0\\
0 & \beta_k
\end{array}\right)
\left(
\begin{array}{cc}
\cdot   & -s\di_y \eps_k(\tau_k) \\
0    &  \cdot
\end{array}\right).
$$ 
This implies
$$
t_k=\alpha_k \cdot s \cdot \di_y \eps_k(\tau_k), 
$$
where $s\asymp -1$.
Now equation 
(\ref{chain rule}) and Lemma~\ref{smalls}(1)  imply
$$
t_k\asymp -\di_y \eps_k (\tau_k) = - \Jac F_n(\tau_k) \asymp   -b^{2^k}  .
$$
\end{proof}

\section{Affine rescaling and quadratic change of variable}\label{quadratic change of variable}

The renormalization procedure described in the previous sections
differs in two ways
from the standard unimodal period-doubling renormalization. 
First, we are renormalizing around the tip of the H\'enon map which becomes the
critical value in the degenerate case. Secondly, we use
non-linear changes of coordinates  $\Psi^n_0$ to define $R^n F$. This was
necessary for the renormalizations to be H\'enon-like maps again.
In this section we will show that in fact, 
a quadratic change of coordinates can be used to produce 
renormalizations converging  to a degenerate universal map. 
(However, affine rescalings would not be sufficient!)
%
This universal map is not the usual fixed point of
 renormalization around the critical point, but rather the fixed point of 
renormalization around the critical value.

Let us now introduce the promised quadratic change of coordinates.
 Take an infinitely renormalizable
 $F\in \II_\Om(\bar\eps)$ with 
sufficiently small $\bar\eps$, so that  the results from
\S~\ref{univ} apply to the maps $\Psi_0^n$. As in that section, 
let us consider translations $T_n: z\mapsto z-\tau_n$ 
(where $\tau_n$ is the tip of $F_n \equiv R^n F$),
and the affine local charts 
$$
L_0^n=  (D^n_0)^{-1} \circ T_0:B_v^{n}(F)\to \mathbb{R}^2.
$$  
Let us represent the maps $F_n$ and $\Psi_0^n$ 
in these charts:   
$$
  \BF_n = T^n\circ F_n \circ (T^n)^{-1}, \quad \BPsi_{0}^{n}=  \id + 
\BS_0^n = L_0^n\circ  \Psi_0^n\circ T_n^{-1}.
$$
%
Let us define the $n^{th}-${\it affine} renormalization of $F$ as follows:
$$
R_{\mathrm{aff}}^n F= L^n_0\circ
 [ F|\, B^{n}_v(F) ]^{2^{n}} \circ (L^n_0)^{-1}= \BPsi_0^n \circ \BF_n \circ (\BPsi_0^n)^{-1}.
$$
Note that the domain of the $n^{th}$-affine renormalizations is the
$\Im \BPsi_0^n$.%
\footnote{Note that $R^n_{\mathrm {aff}}$ is {\it not} the $n$-fold iterate of some $R_{\mathrm {aff}}$.}

 We also let  $T_\infty:z \mapsto z-1$ and 
$$
   \BF_*=  T_\infty \circ  F_*\circ T^{-1}_\infty 
$$ 
By Proposition~\ref{limit}, the maps $ \BPsi^n_0$ converge  to  
$$
       {\bf V}_{*,a_F}: (x,y)\mapsto (\Bv_*(x)+a_Fy^2 , y),
$$ 
exponentially fast. Furthermore, by Theorem~\ref{convergence}, 
$\BF_n\to \BF_*$ exponentially fast. Hence

\begin{thm} \label{affine} 
Let $F\in \II_\Om(\bar\eps)$ be infinitely renormalizable
 with sufficiently small $\bar\eps$. Then
$$
R_{\mathrm{aff}}^n F\to \BV_{*,a_F}\circ \BF_*\circ \BV_{*,a_F}^{-1}
$$
exponentially fast.
\end{thm}

Consider the quadratic change of coordinates $Q_F:\mathbb{R}^2\to \mathbb{R}^2$,
$$
Q_F:(x,y)\mapsto (x-a_Fy^2,y),
$$
and define $H_n: B^{n}_v(F)\to \mathbb{R}^2$ as the composition: 
$$
H_n=Q_F\circ L_0^n.
$$
Conjugating $F^{2^n}$ by these quadratic changes of variable, we obtain the desired renormalizations:
$$
R_{\mathrm{qd}}^n F=H_n\circ F^{2^n}\circ H_n^{-1}.
$$

\begin{figure}[htbp]
\begin{center}
\psfrag{Hn}[c][c] [0.7] [0] {\Large $H_n$}
\psfrag{F}[c][c] [0.7] [0] {\Large $F$}
\psfrag{RnF}[c][c] [0.7] [0] {\Large $R^nF$}
\psfrag{psin}[c][c] [0.7] [0] {\Large $\Psi_0^n$}
\psfrag{T0}[c][c] [0.7] [0] {\Large $T_0$}
\psfrag{Ln}[c][c] [0.7] [0] {\Large $L_0^n$}
\psfrag{Dn}[c][c] [0.7] [0] {\Large $D_0^n$}
\psfrag{Tn}[c][c] [0.7] [0] {\Large $T_n$}
\psfrag{PSIn}[c][c] [0.7] [0] {\Large $\BPsi_0^n$}
\psfrag{FFn}[c][c] [0.7] [0] {\Large $\BF_n$}
\psfrag{raffF}[c][c] [0.7] [0] {\Large $R_{\mathrm{aff}}^nF$}
\psfrag{QF}[c][c] [0.7] [0] {\Large $Q_F$}
\psfrag{RquadF}[c][c] [0.7] [0] {\Large $R_{\mathrm{qd}}^nF$}
\psfrag{Bvn}[c][c] [0.7] [0] {\Large $B_v^n(F)$}
\pichere{0.8}{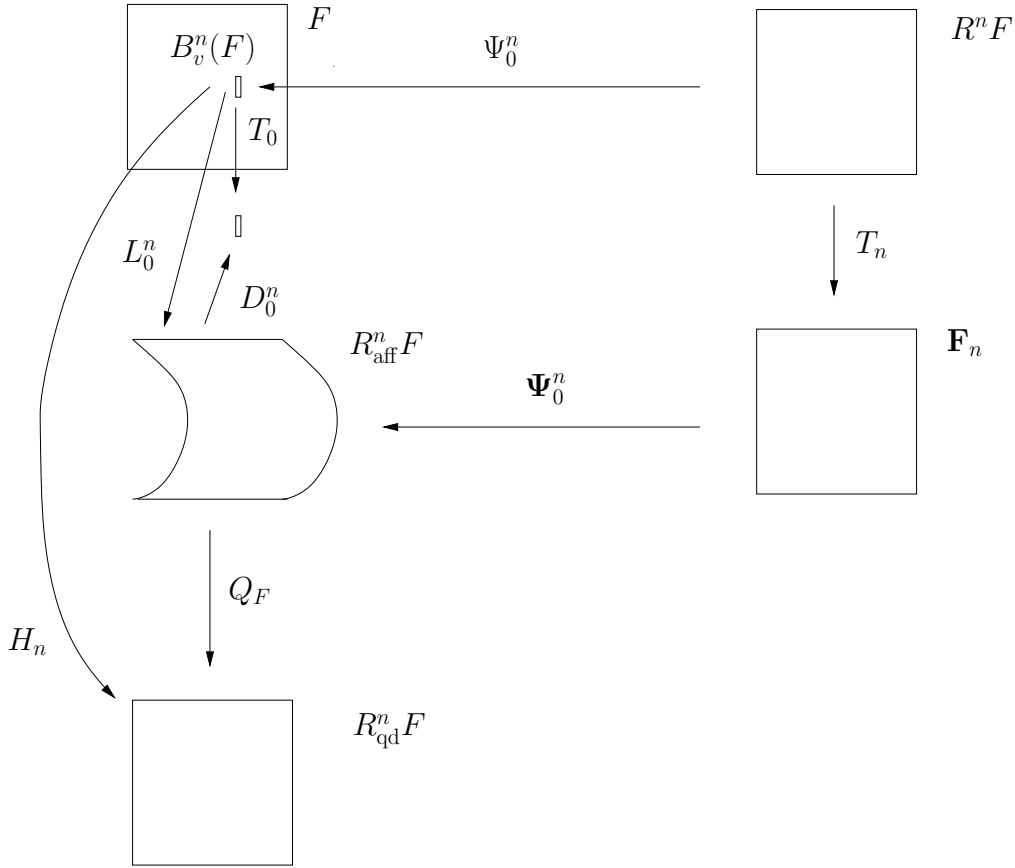}
\caption{Changes of coordinates}
\label{change of coordinates}
\end{center}
\end{figure}

\begin{thm} \label{quad} 
Let $F\in \II_\Om(\bar\eps)$ be infinitely renormalizable
 with sufficiently small $\bar\eps$. Then
$$
R_{\mathrm{qd}}^n F(x,y)\to (l \circ f^*\circ l^{-1}(x), \Bv_*^{-1}(x))
$$
exponentially fast, where $l(x)=(x-1)/ u_*'(1)$.
\end{thm}

\begin{proof} Let $\BV_*=\BV_{*,0}$.  Proposition~\ref{limit} tells us that 
$$
Q_F\circ \BPsi_0^n\to \BV_*,
$$
exponentially fast. This implies that
$$
R_{\mathrm qd }^n = (Q_F\circ \BPsi_0^n) \circ \BF_n \circ (Q_F\circ \BPsi_0^n)^{-1} \to 
\BV_*\circ \BF_*\circ \BV_*^{-1}
$$
exponentially fast. Applying Lemma~\ref{ustarfstar}(3) and the relation between
$u_*$ and ${\bf v}_*$, we obtain:
$$
\BV_*\circ \BF_*\circ \BV_*^{-1}: (x,y)\mapsto (l \circ f^*\circ l^{-1}(x), \Bv_*^{-1}(x)),
$$
and the theorem follows. 
\end{proof}

\begin{rem} In the forthcoming Part II we will construct
 the stable manifold $W^s(\tau_F)$  at the tip $\tau_F$ 
and will show that the number 
$a_F$ is equal to its curvature at $\tau_F$. 
\end{rem}

\begin{rem} The horizontal width of the box $ B^{n}_v(F)$ is proportional 
to the square of its vertical size. This box, a narrow strip containing 
the tip, is aligned along  $W^s(\tau_F)$. 
Any affine change of coordinates   which brings this box roughly to the unit size
is boundedly related to the affine map $L^n_0$.  
In the case when $a_F\ne 0$, these scalings are not capable to ``unbend'' the boxes
$ B^{n}_v(F)$.
(As a model, notice that the rescaling of the parabola $x=ay^2$ by a linear map $(x,y)\mapsto (\si^2 x, \si y)$
does not change the curvature $a$.)
Thus,  the renormalizations obtained by affine changes of variable will always remember
the curvature $a_F$. 
Hence they cannot have a universal limit.   
\end{rem}

\section{Non-existence of continuous invariant line fields}
\label{lineflds}

In this section, $F\in \HH_\Om(\bar\eps)$ stands for an infinitely renormalizable {\it non-degenerate} H\'enon-like map
to which the results of \S \ref{hyp} apply.
Then by the results of \S \ref{Cantor},
it possesses the Cantor attractor $\OO=\OO_F$ on which it acts as the adding machine.
We will show that $F$
does not have continuous invariant line fields on $\OO$. This has several interesting consequences:

\ssk\nin $\bullet$ 
 Contrary to a common intuition,  the attractor $\OO$ does not lie on a smooth curve. 

\ssk\nin $\bullet$
The ${\mathrm {SL}}(2,\R)$-cocycle
\begin{equation}\label{cocycle}
                     z\mapsto DF(z)/ \sqrt{\Jac F(z)}
\end{equation}
is not uniformly hyperbolic over $\OO$.
By Theorem \ref{char exp}, it has non-vanishing characteristic exponents ${\displaystyle \pm \frac{1}{2} \log b}$,
so it is non-uniformly hyperbolic.
 It seems to be the first example of a non-uniformly hyperbolic 
${\mathrm {SL}}(2,\R)$-cocycle over the adding machine.

\begin{lem} \label{lin} If $F$ 
 has a continuous invariant line field on $\OO_F$ then 
there exists $n_0\ge 1$ such that for any $n\ge n_0$, the renormalization $R^nF$ has a 
continuous invariant direction field on $\OO_{R^nF}$. 
\end{lem}

\begin{proof}
 Note first that  any continuous $F$-invariant line field on $\OO=\OO_F$
can be pulled back to a continuous invariant line field on $\OO_{R^n F}$ for any renormalizations $R^n F$.

Furthermore,
since the set $\OO$ is totally discontinuous, 
any continuous invariant line field on it can be  continuously orientated. 
Then there exist a partition of $\OO$ into two 
clopen sets $\OO^+$ and $\OO^-$ such that $F|\OO^+$ preserves the orientation of the field, 
while $F|\OO^-$  reverses it.
Since the pieces $B^n_w$ uniformly shrink as $n\to \infty$,
for  $n$ large enough each $B^n_w\cap \OO_F$  is contained either in $\OO^+$ or in $\OO^-$.  
Hence $F^{2^n}|B^n_v$ either preserves or reverses the orientation of the line field. 
It follows that the renormaliztion $R^nF$ either preserves or reverses the induced 
orientation of the line field on $\OO_{R^nF}$. 
In either case we conclude that the next renormalization, $R^{n+1}F$,  preserves  the induced orientation.
\end{proof}

For any matrix
\begin{equation}\label{A}
   A = \left(
\begin{array}{cc}
a & -\de \\
1 & 0
\end{array}\right), \quad \de>0,
\end{equation}
let us consider its induced action on the circle $S^1$ of directions in $\R^2$.
parametrized by the angle $\theta$. 
(We will keep the same notation, $A$, for the induced action.) 
Let $L$ and $R$ stand for the left- and right-hand semi-circles of $S^1$, while $U$ and $D$
stand for the upper and lower semi-circles. Then
$ A(R)=U,\quad A(L)= D$,
and in the projective coordinate $t=x/y=\ctg \theta$ both maps,  $A: R\ra U$ and $A:L\ra D$, assume the form
\begin{equation}\label{t}
     t\mapsto   a- \frac{\de}{t}.  
\end{equation}

For $\alpha\in (0, \pi/2)$, let us consider two symmetric direction cones: 
$$
    C_\alpha^+ = (\alpha, \pi-\alpha)\equiv  \{ \theta\in S^1: \, \alpha \leq  \theta \leq \pi-\alpha\} ; \quad  
    C_\alpha^- = - C_\alpha^+.      
$$

\begin{lem} \label{V}
There exists an angle $\alpha\in (0, \pi/2)$ with the following property. 
Let $X= \{ F^n(z_0)\}_{n=-\infty}^\infty $ be any two-sided orbit of $F$ in $\OO$, 
and let $z \mapsto \theta(z)$ be an invariant direction field over $X$.
Then there exist points $z^{\pm}\in X$ such that $z^{\pm}\in C^{\pm}_\alpha$.
\end{lem}

\begin{proof} 
Let us  write the differential of $F$ in form (\ref{A}):
$$
   A_z\equiv   DF(z)= \left(
\begin{array}{cc}
a(z) & -\de(z) \\
1 & 0
\end{array}\right).  
$$
Let  $\bar a= \max_{z\in \OO} |a(z)|$  
Without loss of generality  we can assume that $\de(z) < \bar a$ everywhere
(replacing $F$ by its renormalization if needed). 
Let 
$$
   \kappa =  \max \{ 2|\bar a|,\,  1 \}; \quad \alpha =\arcctg \kappa \in (0, \pi/4 ] .
$$ 
We let $Q_i$, $i=1,\dots, 4$, be the four quadrants in $S^1$: 
$$Q_1=[0,\pi/2], \dots,  Q_4= [3\pi/2, 2\pi].$$  

Assume that $\theta(z)\not\in C^+_\alpha$ for any $z\in X$.

Note that  $A_z [0,\alpha] \subset C_\alpha^+ $ for any $z\in \OO$.
Indeed, in the projective coordinate $t=\ctg\theta$, the cone $C^+_\alpha$ is given by  equation $ |t|\leq \kappa$.
By (\ref{t}), we have:   $|\ctg A_z(0)| = |a(z)| < \kappa$ so that $A_z(0)\in C^+_\alpha$.  
If $A_z(\alpha)<\pi/2$, then obviously  $A_z(\alpha)\in C^+_\alpha$ as well.
Otherwise by (\ref{t}) we have:  
$$
     |\ctg A_z(\alpha)| \leq  |\ctg A_z(\pi/4)| = |a(z)-\de(z)| \leq 2 |\bar a| \leq \kappa,
$$
and thus $A_z(\alpha)\in C^+_\alpha$ again.

By invariance of the direction field, we conclude that $\theta(z)\not\in [0,\alpha]$ for $z\in X$.
Hence $\theta(z)\not \in Q_1$ for $z\in X$.     

Since $A_z(Q_4)= [0, A_z(0)]\subset C^+_\alpha\cup Q_1$, we conclude that $\theta(z)\not\in Q_4$.

At this point we already know that $\theta(z)\in [\pi-\alpha,\, \pi]\cup Q_3 \equiv P$ for $z\in X$. 
But  then 
$$ 
    \theta(z) = A_{F^{-1}z} (\theta(F^{-1} z))\subset D,\quad z\in X,
$$   
and hence $\theta(z)\in P\cap D = Q_3$.

By replacing $F$ with its renormalizaion, we can bring  it arbitrary closely to the degenerate fixed point $F_*$.  
Thus, we can assume that the Cantor attractor $\OO_F$ is close to $\OO_{F_*}$ in the first place,
which implies (together with minimality of $\OO_F)$
that $a(z)=f'(z)-\di_x \eps  (z)<0$ for some $z\in X$. But then $A_z (Q_3)\subset Q_4$ for this point $z$, 
and we arrive at a contradiction.  

\ssk We have proved the assertion for the positive cone $C^+_\alpha$.
The one for the negative cone follows by central symmetry of the cocycle.
\end{proof}

\begin{prop} \label{direc} 
There are no continuous invariant direction fields on $\OO_F$.  
\end{prop}

\begin{proof} Suppose there exists a continuous invariant direction field 
on $\OO_F$.  Then there exists such a field for every renormalization. 
By Lemma \ref{V}, for each $n$ we can find a pair of points
$z_n, \zeta_n\in \OO_{R^n F}$ such that
$ \theta(z_n)\in C_\alpha^+ $
while
$
 \theta(\zeta_n) \in C_\alpha^-. 
$

Now project these points to the box $B_v^{n}$ by the maps $\Psi^n_v$
making use of equation (\ref{reshuffling}) and  Lemma~\ref {APPsi}.
We  obtain two sequences of points,  $\hat z_n$ and  $\hat \zeta_n$, converging to  the tip $\tau_F$. 
The direction field at $\hat z_n$  points upward at  angle $\theta(z_n)= \pi/2+ O(b_F)$
while the  direction field at $\hat \zeta_n$  points downward at angle $\theta(\zeta_n)= -\pi/2 +O(b_F)$.
Thus, the direction field is not continuous at the tip of $F$.
\end{proof}

Lemma \ref{lin} and  Proposition \ref{direc}  imply the desired: 

\begin{cor}\label{no line fields}
  The map $F$ does not have a continuous invariant line field on the critical Cantor set $\OO_F$.
\end{cor} 

It immediately yields:

\begin{thm} \label{nonunivhyp} The map $F$  is not partially  hyperbolic on $\OO_F$
in the sense that the contracting and neutral line fields corresponding to the characteristic
exponents $\log b$ and $0$ (see Theorem \ref{char exp}) are discontinuous.
\end{thm}

\begin{thm}\label{cocycleTh}  The ${\mathrm {SL}}(2,\R)$-cocycle (\ref{cocycle})
is non-uniformly hyperbolic over $\OO$.
\end{thm}

\begin{thm} \label{curve}
There are no smooth curves containing $\OO_F$.
\end{thm}

\begin{proof}
 If $\CC$ is a smooth curve containing $\OO_F$, then its tangent lines $l(z)$
give us a continuous line field on $\OO_F$.
 Since $\OO_F$ does not have isolated points,  
$$
 l(z)= \lim_{\zeta\to z} l(z,\zeta),  
$$
where $l(z,\zeta)$ is the line passing through $z$ and $\zeta\in \OO_F$, $\zeta\not=z$.
It follows that the line field $l(z)$ is invariant over $\OO_F$, contradicting Corollary~\ref{no line fields}. 
\end{proof}

\section{Non-rigidity of the critical Cantor set}\label{non-rigidity}

We will show that the invariant Cantor set $\OO$ of an infinitely
renormalizable H\'enon-like map is not rigid. In fact,
there is a definite upper bound smaller than 1 on the H\"older exponent
of the conjugacy
between two such Cantor sets of any two H\'enon-like maps with different average
Jacobians. 

\begin{thm}\label{opthol}  Let $F$ and $ {\tilde F}$ be two infinitely renormalizable
H\'enon-like maps with average Jacobian $b$ and ${\tilde  b}$ resp.
Assume $b> {\tilde b}$. Let $\phi$ be a homeomorphism which conjugates $
F|_{\OO_F}$ and ${\tilde F}|_{\OO_{{\tilde F}}}$ with $\phi(\tau({\tilde F}))=\tau(F)$. Then
the H\"older exponent of $\phi$ is at most $\frac12(1+ \ln b / \ln \tilde b )$.
\end{thm}

\begin{proof} We let $F_k=R^kF$ be the $k$-fold renormalization,  
$v_k=\tau(F_k)$ be its tip, 
$c_k=(F_k)^{-1}(v_k)$ be its ``critical point'',  
and $c_{k}^{k+n}=\Psi^{k+n}_k(c_{k+n})$. 
Furthermore, let $w_k=F_k (v_k)$ and $z_{k}^{k+n}=F_k (c_{k}^{k+n})$,
see Figure~\ref{nnrg}. 
We will mark the corresponding objects of ${\tilde F}$ with the tilde.

\begin{figure}[htbp]
\begin{center}
\psfrag{RknF}[c][c] [0.7] [0] {\large $R^{k+n}F$}
\psfrag{F}[c][c] [0.7] [0] {\large $F$}
\psfrag{RkF}[c][c] [0.7] [0] {\large $R^kF$}
\psfrag{Pk}[c][c] [0.7] [0] {\large $\Psi_0^k$}
\psfrag{Pkn}[c][c] [0.7] [0] {\large $\Psi_k^{k+n}$}
\psfrag{ckn}[c][c] [0.7] [0] {\large $c_{k+n}$}
\psfrag{vkn}[c][c] [0.7] [0] {\large $v_k^{k+n}$}
\psfrag{vk}[c][c] [0.7] [0] {\large $v_k$}
\psfrag{vkn}[c][c] [0.7] [0] {\large $v_{k+n}$}
\psfrag{wk}[c][c] [0.7] [0] {\large $w_k$}
\psfrag{zkn}[c][c] [0.7] [0] {\large $z_k^{k+n}$}
\psfrag{Zkn}[c][c] [0.7] [0] {\large $Z_k^{k+n}$}
\psfrag{Wk}[c][c] [0.7] [0] {\large $W_k$}
\psfrag{ckkn}[c][c] [0.7] [0] {\large $c_k^{k+n}$}
\pichere{1.0}{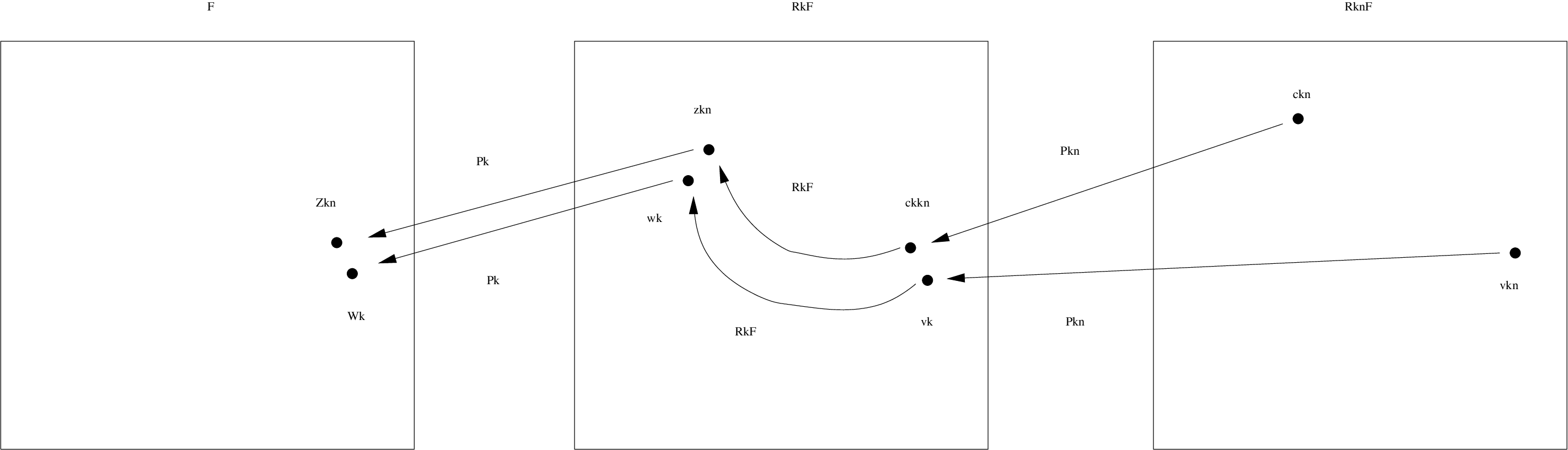}
\caption{}
\label{nnrg}
\end{center}
\end{figure}

For large  renormalization levels $k\ge 1$, we have:  $b^{2^k}\gg {{\tilde b}}^{2^k}$.
Choose now the scale $n=n(k) \ge 1$ satisfying
$$
      \si^{n+1} \leq  {\tilde b}^{2^k} < \si^n.
$$
Let $\Delta  {\tilde x}$ and $\Delta  {\tilde y}$ be the differences between the
$x$- and $y$-coordinate of the points $ {\tilde v}_k$ and ${\tilde c}_{k}^{k+n}$.
Representation (\ref{reshuffling}),  Lemma~\ref{APPsi} and Corollary~\ref{tilt} imply:
$$
|\Delta {\tilde y}| \asymp \si^n
$$
and
$$
|\Delta {\tilde x}|=O( \si^{2n}+{{\tilde b}}^{2^k}\cdot 
|\Delta {\tilde y}|)= O( \si^{2n}).
$$
Applying ${\tl F}_k$ to these points using the Universality  Theorem~\ref{universality}, we obtain :
$$
\begin{aligned}
\dist({\tilde z}_{k}^{k+n}, {\tilde w}_k)&= O( |\Delta {\tilde x}|+ 
|\Delta {\tilde y}|\cdot \frac{\partial \eps_k}{\partial y})\\
&
= O(\si^{2n} + \si^n {\tilde b}^{2^k})  =O(\si^{2n} ). 
\end{aligned}
$$
(Notice that ${\tl F}_k$ has compressed the vertical distance between $ {\tilde v}_k$ and ${\tilde c}_{k}^{k+n}$
to make  it comparable with the horizontal distance.)

Consider now points ${\tilde Z}_{k}^{k+n}=\Psi_0^{k}({\tilde z}_{k}^{k+n})$ and
${\tilde W}_k=\Psi_0^{k}({\tilde w}_k)$ in the domain of ${\tilde F}$. 
By Lemma 5.1, we have:  
$$
\dist({\tilde W}_k, {\tilde Z}_{k}^{k+n})=O(\si^{2n+k}).
$$ 

Let us now estimate the distance between the corresponding points
for $F$. For the same reason as above, we have: $ |\Delta y|\asymp \si^n.$
Furthermore, since the tilt of the box $B^{n+k}_k$ is of order $b^{2^k}$ (by Corollary~\ref{tilt}), 
we obtain for some $\gamma>0$:
$$
|\Delta x|\ge 2 \gamma  \left( b^{2^k} |\Delta  y|- \si^{2n} \right) \geq \gamma\, b^{2^k} \si^n,
$$
where the last estimate uses that $b^{2^k}\gg \si^n$. 
Hence
$$
  | \pi_2(w_k) - \pi_2(z_{k}^{k+n})| = |\Delta x|\ge \gamma \,  b^{2^k}\si^n, 
$$
where $\pi_2$ stands for the vertical projection.
Using  representation (\ref{reshuffling}) and  Lemma~\ref{APPsi} once again, we obtain: 
$$
\dist(W_k, Z_{k}^{k+n}    ) \ge \gamma\, \si^{k+n} b^{2^k}
$$

Any H\"older exponent $\alpha>0$ for the conjugating
homeomorphism has to satisfy
$$
\dist(W_k, Z_{k}^{k+n})\le C \,
 (\dist({\tilde W}_k, {\tilde Z}_{k}^{n} ))^\alpha.
$$
Hence
$$
\si^k\,  { {\tilde b}}^{2^k} \, b^{2^k} \le
C\, \left( \si^k \, {{\tilde b}}^{2^k}\, {\tilde b}^{2^k}  \right)^\alpha
$$
which implies the asserted bound:
$$
\alpha\le \frac12 \left(1+\frac{\ln b}{\ln  {\tilde b}}\right).
$$
\end{proof}

\begin{cor}  Let $F$ be an infinitely renormalizable H\'enon-like
map with the average Jacobian $b$ and $F_0$ be a degenerate
infinitely renormalizable H\'enon-like map.
Let $\phi$ be a homeomorphism which conjugates $ F|_{\OO_F}$ and $
F_0 |_{\OO_{F_0}}$ with $\phi(\tau( F_0))=\tau(F)$. Then the H\"older
exponent of $\phi$ is at most $\frac12$. \qed
\end{cor}

\section{Generic unbounded geometry}\label{unbounded geometry sec}

 An infinitely renormalizable H\'enon map has {\it bounded geometry} if
$$
{\diam}(B^n_{w\nu})\asymp \dist(B^n_{w v},B^n_{w c}),
$$
for $n\ge 1$ and $w\in W^{n-1}$ and $\nu\in W$. A
slight modified version of this definition would require
$$
{\diam}(B^n_{w \nu}\cap \OO)\asymp \dist(B^n_{w
v}\cap \OO,B^n_{w c}\cap \OO).
$$
The following theorem holds for both definitions, with the same proof:

\begin{thm}\label{unbdgeomth} 
Let $F_b$, $b\in [0,1]$, be a family of infinitely
renormalizable H\'enon-like maps parameterized by the average Jacobian,
that is,  $b_{F_b}=b$. Then for some $b_0>0$, the set of parameter values for which
$F_b$ does not have bounded geometry contains a dense $G_{\delta}$ subset in 
an 
interval $[0,b_0]$.
\end{thm}

\begin{proof} Let us take $\bar b>0$ so small that the estimates of \S 7 
on $\Psi_k^n$ hold for all $F_b$ with $b\in [0, 2\bar b]$. For $n >  k\geq 1$, 
let us consider the boxes  $B_{k}^{n}=\Psi_{k}^{n}(B)$ in the
domain of $F_k\equiv R^kF$, and let   
$$
P_{k}^{n}=\Psi_{k}^{n-1} (F_{n-1}(B_{n-1}^{n})).%
\footnote{In notations of \S \ref{pieces}, $B^n_k= B^{n-k}_{v^{n-k}}(F_k)$, $P^n_k= B^{n-k}_{v^{n-k-1}c}(F_k)$.} 
$$
 Note that $ B_{k}^{n}\cup P_{k}^{n}\subset B_{k}^{n-1}.$
 As in \S \ref{Psi-functions}, $\tau_k=\tau(F_k)$ stands for the tip of $F_k$. 
 Let us also consider some point $c_k\in P_k$ moving continuously with the parameter
 (for instance, we can take the ``critical point'' $c_k=(F_k)^{-1}(v_k)$ of $F_k$),  
and let  $c_{k}^{n}=\Psi^{n}_k(c_{n})\in B^n_k$
(compare Figure 10.1). 

Making use of representations  (\ref{reshuffling}) and (\ref{factoring}),
let us estimate the relative horizontal positions of the points $\tau_k$ and $c^n_k$.
Let 
$$
   z=(x, y)= (\id + {\bf S}^n_k) (\tau_n), \quad z_0=(x_0,y_0) = (\id + {\bf S}^n_k) (c_n).
$$ 
By Lemma \ref{ustar}, we have:
$$
   x - x_0 =  {\bf v}_*(c_n)-{\bf v}_*(\tau_n) + O(\bar b^{2^k}+\rho^{n-k}),   
$$
which is a negative number of order 1, provided $k$ and $n-k$ are sufficiently big ($\geq N$).
Hence
$$
  \pi_1(c^n_k)- \pi_1(\tau_k) = \pi_1(D^n_k (z-z_0))
$$
$$
  = [ \si^{2(n-k)} (x-x_0) + t_k (-\si)^{n-k} (y-y_0)]\, (1+O(\rho^k))
$$
Together with Corollary~\ref{tilt}, the above estimates yield for even $n-k$:
\begin{equation}\label{hor proj} 
    \pi_1( c^n_k) - \pi_1(\tau_k)  = \si^{2(n-k)} (x-x_0)[ 1  - b^{2^k}\si^{-(n-k)} r_{n,k}]\, (1+O(\rho^k),  
\end{equation}
where  $0< r \leq r_{n,k}\le \rho $ uniformly in $b$.  

Let us now take any parameter $b_- \in (0, \bar b)$ and any integer  $k\geq N$. 
Let us find the biggest  $n$ such that $n-k$ is even  and $\si^{n-k}> \rho (b_-)^{2^k}$.
By (\ref{hor proj}), for the map $F_{b_-}$, the point $c^n_k$ lies to the left of the tip $\tau_k$. 
Let us increase $b_-$ to a parameter $b_+$ such that $(b_+)^{2^k} = 2 r^{-1} \si^{n-k}$.
Then for $F_{b_+}$, the point $c^n_k$ lies to the right of the tip $\tau_k$.
Hence there exists a parameter $b\in (b_-, b_+)$  for which $c_n^k$ lies strictly below the tip $\tau_k$.   

Moreover,
\begin{equation}\label{b vs si}
          b^{2^k}\asymp \si^{n-k},
\end{equation}
and the hyperbolic distance between $b$ and $b_-$ in the hyperbolic line $\R_+$ is small:  
$\ln (b / b_-)=O(2^{-k})$. Letting $k$ run through all integers $N, N+1, \dots$,
we obtain a dense set of parameters $b\in (0 , \bar b)$ for which the point $c^n_k$ lies strictly below
the tip $\tau_k$ for some $k, n$. It follows that there is a open and dense subset $\La_k\subset (0, \bar b)$
of parameters for which  some point $c^n_k\in P^n_k$%
\footnote{We keep the same notation for this point, though it is not necessarily the one chosen above}
lies strictly below the tip $\tau_k$
for some $n>k$. Hence for any parameter $b$ in the open $ G_\de$-set $\La= \cap\La_k$, 
this happens for infinitely many levels $k$.

 \msk
We are going to show that the geometry of the critical Cantor set degenerates for $b\in \La$. 
It is convenient to shift the level by 1, so that we assume that $b\in \La_{k+1}$. 
Let $w_k$ and $z^n_k$  be  the images of the points $\tau_{k+1}$ and $c^n_{k+1}$
under the the map $F_k\circ \Psi^{k+1}_k$ (which is equal to $\Psi_c^1(F_{k+1})$ in notation of \S \ref{branches}).
Since the maps $\Psi_c^1$ preserve the vertical foliation (see Remark \ref{vertical fol preserved}),
the points $w_k$ and $z^n_k$ also lie one strictly  above the other.

Since the point $c^n_{k+1}$ lies strictly below $\tau_{k+1}$ on distance of order $\si^{n-k}$, 
the interval between the points $\Psi^{k+1}_k (c^n_{k+1})$ and $\Psi^{k+1}_k(\tau_{k+1})$ has length of order $\si^{n-k}$
and slope of order $-b^{2^k}$ (see  Lemma \ref{smalls}). Hence the distance between the horizontal projections of these two points
is of order $\si^{n-k}b^{2^k}$. But it is equal to the distance between their $F_k$-images, $z^n_k$ and $w_k$.
Thus, 
$$
    \dist(w_k, z^n_k)\asymp \si^{n-k} b^{2^k}.
$$
  
Applying $F_k$ once more, we obtain two point on the same horizontal line such  that
\begin{equation}\label{dist}
\dist(F_k(w_k), F_k(z^n_k))\asymp \sigma^{n-k} \,  b^{2^{k+1}}.
\end{equation}

\ssk
Let us now estimate the sizes of the corresponding pieces.
Let $Q$ stand for either $B^n_k$ or $P^n_k$. 
By  (\ref{reshuffling}), Proposition~\ref {limit} and Corollary \ref{tilt},
it contains two points such that the interval joining them
has length of order $\si^{n-k}$ and tilt of order $b^{2^k}$.
Hence 
$$
 |\pi_1(Q)|\geq \gamma \si^{n-k} b^{2^k}
$$
for some $\gamma>0$.
It follows that both projections of $F_k^2(Q)$ are at least that big (up to a constant). 
We are interested only in the vertical size:
$$
  | \pi_2(F_k^2 (Q))|\geq  \gamma \si^{n-k} b^{2^k}.
$$
Comparing this with (\ref{dist}), we see that the distance between the points
$F_k(w_k)$ and  $F_k(z^n_k)$ 
is at least $b^{2^k}$ times smaller than the vertical size of the pieces
$F_k^2(B^n_k)$ and $F_k^2(P_n^k)$ that contain these points.

\ssk
Finally, we should bring these two pieces to the domain of $F$
by the map $\Psi^k_0$. Since this map contracts the horizontal distances
stronger than the vertical ones, the gap between the images of the pieces
will be even smaller compared to the size of the pieces
(the gap will become at least $b^{2^k} \si^k$ times smaller than the size of the pieces). 

The conclusion follows.
\end{proof}

\section{H\"older geometry of the critical Cantor set}\label{holder}

If $P=B^{n-1}_{\sigma}$, $n\ge 1$ and $\sigma\in \Sigma^{n-1}$, is a piece of an infinitely 
renormalizable H\'enon-like map  
$F\in \II_\Om(\bar\eps)$ we call
the distance $g=\dist(B^n_{\sigma v},B^n_{\sigma c})$ the 
{\it gap} of the piece $P$. An infinitely renormalizable H\'enon map has {\it H\"older  bounded geometry} if
there exist $\alpha>0$  and $C>0$ such that
$$
g^\alpha \ge C \cdot {\diam}(P)
$$
for every piece $P$ of $F$.
\begin{thm}\label{holgeo} Every infinitely renormalizable H\'enon-like map  $F\in \II_\Om(\bar\eps)$, with 
sufficiently small $\bar\eps$, has H\"older bounded geometry.
\end{thm}

The proof of this Theorem will be by induction in the size of the pieces. The beginning of the induction is the 
following Proposition.

\begin{prop} \label{Ind0} There exist constants $K, C>0$ such that for every  $F\in \II_\Om(\bar\eps)$ with 
sufficiently small $\bar\eps$ and every piece $P$ of $F$ with gap $g$ the following holds. If
$$
\mathrm{diam}(P)\ge K \cdot b_F
$$
then
$$
g\ge C\cdot \mathrm{diam}( P) .
$$
\end{prop}

\begin{rem} In the previous section we showed that the geometry of $\OO_F$
might be 
unbounded. Proposition ~\ref{Ind0} states that this two-dimensional phenomenon 
becomes observable only at a scale of the order of $b$. 
\end{rem}

The proof of this Proposition relies on the following Lemma for which we need some notation. 
Given a piece $P$, let $H$ and $V$ stand for its 
horizontal and vertical projections. Let
$$
q_P=\frac{|V|}{|H|}.
$$
The piece $P$ is obtained by repeatedly applying contractions, 
say $P=\Psi_{\omega_1\omega_2\dots \omega_{n}}^n(B)$. Let 
$P_k=  \Psi_{\omega_k\omega_{k+1}\dots \omega_{n}}^n(B)$ be the 
corresponding piece of $F_k\equiv R^kF$, $k\le n$.

\begin{lem}\label{q} For every $K>0$ there exists $C>0$ such that if $P$ is a piece of $F$ with
$$
\mathrm{diam}(P)\ge K \cdot b_F
$$
then
$$
q_k=q_{P_k}\le C\cdot \frac{1}{b_F},
$$
for $k\ge 1$.
\end{lem}

\begin{proof} The piece $P$ is of the $n^{th}$
generation of $F$.  Let $1\le k\le n$ and $s\ge k$ be maximal such that
$$
P_k=\Psi^{s-k}_{v^{s-k}}(P_s),
$$
(where only ``critical value'' contractions were used). Then
$$
P_s=\Psi_{c\omega_{s+1}\dots \omega_{n}}^{n-s}(B)
$$
Let 
$$
P'=\Psi_{v\omega_{s+1}\dots \omega_{n}}^{n-s}(B)\subset B^1_v(F_s).
$$
Note,
$$
F_s(P')=P_s.
$$
Let $H_s, V_s$ and $H', V'$ be the horizontal and vertical projections of $P_s$ and $P'$ respectively.
From Theorem~\ref{universality}, for some uniform $A>0$ and $K_1>0$ 
$$
K \cdot b\le \mathrm{diam}(P_s)\le |V_s| + |H_s| \le |V_s|+ A |H'|+K_1b^{2^s}.
$$
Because $|V_s|=|H'|$ we get
\begin{equation}\label{H'}
|H'|\ge K_3 \cdot b,
\end{equation}
for some $K_3>0$. From Theorem~\ref{universality} we get for some uniform 
$a>0$ and $K_4>0$ 
\begin{equation}\label{Hs}
|H_s|\ge a |H'|-K_4b^{2^s}.
\end{equation}
Now ~\ref{H'} and ~\ref{Hs} imply
$$
q_s=\frac{|V_s|}{|H_s|}\le \frac{|H'|}{a|H'|-K_4 b^{2^s}}=O(1).
$$
From  Proposition~\ref {limit} and (\ref{reshuffling}) 
we get
$$
q_k= O(1/\sigma^{s-k}).
$$
Using Lemma~\ref{contracting} 
$$
\begin{aligned}
b_F&\le \frac1K\cdot \mathrm{diam}(P)\le \frac1K\cdot \mathrm{diam}(P_k)\\
&\le \frac1K\cdot\mathrm{diam}(P_s) \cdot C\sigma^{s-k} \le
 \frac{C}{K} \cdot  \sigma^{s-k}.
\end{aligned}
$$
And the Lemma follows.
\end{proof}

\noindent
{\it Proof  of Proposition~\ref{Ind0}.}
Let $P_k$ be a piece (of some $F_k$) of generation $n-k$. Let $G_h\subset H$ and 
$G_v\subset V$ be the minimal closed
intervals such that $G_h\times V$ and $H\times G_v$ do intersect the two 
pieces of the next generation contained in $P_k$. Note, $G_v$ (and $G_h$) is a 
degenerate interval if the pieces of the next generation have intersecting 
vertical (horizontal) projections. The following argument will show that this
does not happen. Let
$$
\Gamma_{k,n}^{\text{hor}}=\min_{P_k} \frac{|G_h|}{|H|}
$$
$$
\Gamma_{k,n}^{\text{ver}}=\min_{P_k} \frac{|G_v|}{|V|}
$$
and
$$
\Gamma_{k,n}=\min\{ \Gamma_{k,n}^{\text{hor}}, \Gamma_{k,n}^{\text{ver}}
\}.
$$
 Let $\mathcal{P}_{k,n}$ be the pieces of generation $n-k$ of $F_k$ and
$$
\mathcal{P}^c_{k,n}=\{P\in\mathcal{P}_{k,n}| P\in B^1_c(F_k)\}
$$
and
$$
\mathcal{P}^v_{k,n}=\{P\in\mathcal{P}_{k,n}| P\in B^1_v(F_k)\}.
$$
Also define
$$
\Gamma_{k,n}^{\text{hor}, c}=
\min_{P_k\in \mathcal{P}^c_{k,n}} \frac{|G_h|}{|H|},
$$
$$
 \Gamma_{k,n}^{\text{hor}, v}=
\min_{P_k\in \mathcal{P}^v_{k,n}} \frac{|G_h|}{|H|}.
$$
And 
similarly, define $ \Gamma_{k,n}^{\text{ver}, c}$ and
 $\Gamma_{k,n}^{\text{ver}, v}$. Observe, using the specific normalization
of H\'enon-like maps (y'=x) and the fact that the functions $\psi^1_v(F_k)$ 
are affine in the vertical direction,
\begin{enumerate}
\item $\Gamma_{k,n}^{\text{ver}, v}=\Gamma_{k+1,n}^{\text{ver}}$,
\item $\Gamma_{k,n}^{\text{ver}, c}= \Gamma_{k,n}^{\text{hor}, v}$,
\item $\Gamma_{k,n}^{\text{hor}, c} \ge \Gamma_{k,n}^{\text{ver}, v}$.
\end{enumerate}
The last property follows from Lemma ~\ref{permute} (3).
These relations imply
\begin{equation}\label{gkn}
\Gamma_{k,n}\ge \min\{\Gamma_{k+1,n}^{\text{ver}}, 
\Gamma_{k,n}^{\text{hor}, v}\}.
\end{equation}
Now we will express $\Gamma_{k,n}^{\text{hor}, v}$ in terms of 
$\Gamma_{k+1,n}^{\text{hor}}$. 
Let $P\in \mathcal{P}_{k+1,n}$ and
$G_h\subset H$ and $V$ be the corresponding intervals. 
Let $\hat{P}=\psi^{k+1}_v(P)$ and 
$\hat{G}_h\subset \hat{H}$.  Then, using Lemma~\ref{smalls},
(\ref{Dk}), and the tilt quantified in Corollary~\ref{tilt} 
$$
|\hat{G_h}|\ge D_g |G_h|-K_1\cdot |V| \cdot b^{2^{k+1}},
$$
and
$$
|\hat{H}|\le D_h |H|+K_1\cdot |V| \cdot b^{2^{k+1}},
$$
where
$$
D_g=\frac{\partial \Psi_v^{k+1}}{\partial x}(x_g,y_0),
$$
$$
D_h=\frac{\partial \Psi_v^{k+1}}{\partial x}(x_h,y_0),
$$
with $x_g\in G_h$, $x_h\in H$ appropriately chosen, $y_0\in \partial V$, and 
$K_1>0$. Lemma~\ref{smalls}(3) and Lemma~\ref{contracting} gives
$$
\ln \frac{D_g}{D_h}=O(\sigma^{n-k}).
$$
These estimates, together with Lemma \ref{q} and the assumption that $\mathrm{diam}(P)\ge K\cdot b$, 
imply that for some constant $K_2, K_3>0$
$$
\frac{|\hat{G_h}|}{|\hat{H}|}\ge \frac{|G_h|}{|H|}
                                 \cdot \exp(-K_2\cdot \sigma^{n-k})\cdot
                                 \frac{1- K_3\cdot b^{2^{k+1}-1}\cdot 
\frac{|H|}{|G_h|}} {1+ K_3\cdot b^{2^{k+1}-1}}.
$$
This implies
\begin{equation}\label{iterg}
\Gamma_{k,n}^{\text{hor}, v}\ge \frac{e^{-K_2\sigma^{n-k}}}{1+ K_3\cdot b^{2^{k+1}-1}}\cdot
\left[
\Gamma_{k+1,n}^{\text{hor}}- K_3\cdot b^{2^{k+1}-1}
\right].
\end{equation}
Equation (\ref{gkn}) and (\ref{iterg}) imply
\begin{equation}\label{g}
\Gamma_{k,n}\ge \frac{e^{-K_2\sigma^{n-k}}}{1+ K_3\cdot b^{2^{k+1}-1}}\cdot
\left[
\Gamma_{k+1,n}- K_3\cdot b^{2^{k+1}-1}
\right].
\end{equation}
By iterating estimate (\ref{g}) and using that $\Gamma_{n-1,n}\asymp 1$ we get 
$m>0$ such that
$$
\Gamma_{0,n}\ge m >0,
$$
for $n\ge 1$.
This implies Proposition~\ref{Ind0}.
\qed

\bigskip

The induction hypothesis (denoted by $\mathrm{Ind}_n$, $n\ge 0$) we will
 use to prove Theorem \ref{holgeo} is the following.
There exist $\alpha_n>0$ and constants $C>0$ and $K>0$,
 independent of $F$ and  $n\ge 0$, such that the condition
$$
\mathrm{diam}(P)\ge K \cdot b^{2^n},
$$
on any piece $P$ of $F$ implies
$$
g^{\alpha_n}\ge C \cdot \mathrm{diam}(P).
$$ 
Proposition \ref{Ind0} states that $\mathrm{Ind}_0$ holds with $\alpha_0=1$.

\bigskip

Assume that $\mathrm{Ind_j}$ holds for $j\le n$. We are going to prove $\mathrm{Ind}_{n+1}$.
Consider a piece $P_{n+1}$ of $F$  with
$$
\mathrm{diam}(P_{n+1})\ge K \cdot b^{2^{n+1}}.
$$
Because $\mathrm{Ind_j}$ holds for $j\le n$ we may assume without loss of generality that  $\mathrm{diam}(P_{n+1})\le K \cdot b^{2^{n}}$.
This piece is obtained by applying a contraction $\Psi^1_{c}(RF)$ or $\Psi^1_{v}(RF)$ to a piece $P_n$ of $RF$.
Note that
$$
\mathrm{diam}(P_n)\ge \mathrm{diam}(P_{n+1})
\ge K (b^2)^{2^n}.
$$
Hence, if $g_n$ is the gap of $P_n$, $\mathrm{Ind}_n$ implies
$$
g_n^{\alpha_n}\ge C \cdot \mathrm{diam}(P_n).
$$
Observe,
$$
g_{n+1}\ge A\cdot b \cdot g_n,
$$
for some constant $A>0$. 
We need to find an estimate for $\alpha_{n+1}>0$ such that
\begin{equation} \label{cond}
g_{n+1}^{\alpha_{n+1}}\ge C \cdot \mathrm{diam}(P_{n+1}).
\end{equation}
We may assume $\alpha_{n+1}\le \alpha_n$.
The condition ~\ref{cond} holds if
\begin{equation}\label{cond2}
 (A \cdot b)^{\alpha_{n+1}}\cdot (C\cdot\mathrm{diam}(P_n) )^{\frac{\alpha_{n+1}}{\alpha_n}}\ge C\cdot \mathrm{diam}(P_n).
\end{equation}  
Use the fact that for some $L>0$
$$
\mathrm{diam}(P_n)\le L\cdot \frac{1}{b} \cdot \mathrm{diam}(P_{n+1})\le 
\frac{L}{K}\cdot b^{2^{n}-1}
$$
to reduce the  condition ~\ref{cond2} to the next sufficient condition for 
\ref{cond}. 
Namely,
\begin{equation}\label{cond3}
A^{\alpha_{n+1}}\ge (C \cdot L)^{1-\frac{\alpha_{n+1}}{\alpha_n}} \cdot
b^{(2^n-1)\cdot (1-\frac{\alpha_{n+1}}{\alpha_n})-\alpha_{n+1}}.
\end{equation}
Finally, this condition ~\ref{cond3} reduces to the sufficient condition
$$
-M\ge \ln b\cdot [(1-\frac{\alpha_{n+1}}{\alpha_n})\cdot (2^n-1)-1],
$$
where $M>0$ is some large constant.
Now choose $\alpha_{n+1}$ such that
$$
(1-\frac{\alpha_{n+1}}{\alpha_n})\cdot (2^n-1)=m
$$
is constant but sufficiently large and one obtains $\alpha_{n+1}>0$ for which $\mathrm{Ind_{n+1}}$ holds.
Moreover, the sequence $\alpha_n>0$ decreases to some $\alpha>0$. This finishes the proof of the 
Theorem \ref{holgeo}.

\section{Open Problems}\label{problems}

Let us  finish with some further questions that naturally arise from the previous discussion.
The first two of them are probably very hard, 
while others should be more tractable.

\begin{enumerate}

\item Prove that $F_*$ is the only fixed point of the H\'enon renormalization $R$,
   and $R^n F\to F_*$ exponentially  for any infinitely renormalizable H\'enon-like
   map $F$. 

\item Is it true that the trace of the unstable manifold $\WW^u(F_*)$
  by the two-parameter H\'enon family $F_{c,b}: (x,y)\mapsto (x^2+c-by, x)$
  is a (real analytic) curve $\gamma$ on which the Jacobian $b$ assumes all
  values $0<  b<1$. If so, does this curve converge to some particular point 
  $(c,1)$ as $b\to 1$?  

\item How good is the conjugacy $h\colon \OO_F\to \OO_G$ when
$b_F=b_G$?

\item Is the conjugacy $h\colon \OO_F\to \OO_G$ always H\"older?
An equivalent question (due to Theorem \ref{holgeo}) 
is whether the pieces $B^n_\si$ decay no faster than exponentially in $n$? 
The answer is probably negative in general. 

\item 
Can $\OO_F$ have bounded geometry when $b_F\ne 0$?
If so, does this property depend only on the average Jacobian $b_F$?  

\item Does the Hausdorff dimension of $\OO_F$ depend only on the average Jacobian $b_F$?
(This question was suggested by A. Avila.)
\end{enumerate}

\section{Appendix: Shuffling}
\label{shuff}

In this section we will briefly recall  some analysis of long compositions of 
diffeomorphisms of the interval. It is convenient to represent a $C^3$ 
diffeomorphism $\phi:[-1,1]\to [-1,1]$ by its $C^1$ non-linearity
$$
\eta_\phi=\frac{D^2\phi}{D\phi}.
$$ 
The following  Lemma was used in \S\ref{univ}.

\begin{lem}\label{shufflem}(Shuffling)
 For every $B>0$ there exists $K>0$ such that the following holds.
Let $\phi_j, \phi_j^*:[-1,1]\to [-1,1]$, $j=1,\dots, n$ be $C^3$ 
diffeomorphisms and let
$$
\Phi=\phi_n\circ \dots \circ \phi_2\circ \phi_1
$$
and
$$
\Phi^*=\phi^*_n\circ \dots \circ \phi^*_2\circ \phi^*_1.
$$
If
$$
\sum_{j=1}^n \|\eta_j\|_{C^1}\le B 
$$
and 
$$
\sum_{j=1}^n \|\eta^*_j\|_{C^1}\le B
$$
where $\eta^{(*)}_j$ is the non-linearity of $\phi^{(*)}_j$, then
$$
\mathrm{dist}_{C^2}(\Phi, \Phi^*)\le K \sum_{j=1}^n \|\eta_j-\eta^{*}_j\|_{C^0}.
$$
\end{lem} 

This Lemma is a consequence of the Sandwich-Lemma 10.5 from \cite{Ma}.
Here we will use a slightly different version of this Sandwich-Lemma, whose 
proof is exactly the same as the proof for the original formulation.

\begin{lem}(Sandwich)\label{Sandwich} For every $B>0$ there exists $K>0$ such that the 
following holds.
 Let $\phi_j, \phi:[-1,1]\to [-1,1]$, $j=1,\dots, n$ be $C^3$ 
diffeomorphisms and let
$$
\Phi=\phi_n\circ \dots \circ \phi_{k+1}\circ \phi_k\circ \cdots
 \circ \phi_2\circ \phi_1
$$
and
$$
\Psi=\phi_n\circ \dots \circ \phi_{k+1}\circ \phi\circ \phi_k\circ \cdots
 \circ \phi_2\circ \phi_1.
$$
If
$$
\sum_{j=1}^n \|\eta_{\phi_j}\|_{C^1} +\|\eta_{\phi}\|_{C^1}\le B, 
$$
 then
$$
\|\eta_\Phi -\eta_\Psi\|_{C^0}\le K \|\eta_\phi\|_{C^0}.
$$
\end{lem}

The proof for the Shuffling-Lemma \ref{shufflem} 
consists of {\it sandwiching} the 
diffeomorphisms $\phi_k^*\circ \phi_k^{-1}$ between $\phi_{k+1}$ and $\phi_k$,
$k=1,\dots, n$. In this way $\Phi$ is changed into $\Phi^*$. To estimate the distance between these two diffeomorphism we need the following Lemma.

\begin{lem} \label{inv} For every $B>0$ there exists $K>0$ such that the following holds.
 Let $\phi, \psi:[-1,1]\to [-1,1]$ be $C^3$ diffeomorphisms with
$$
\|\eta_\phi\|_{C^0} \le B
$$
Then
$$
\|\eta_{\psi\circ \phi^{-1}}\|_{C^0}\le K\cdot 
\|\eta_\psi-\eta_\phi\|_{C^0}
$$
and 
$$
\|\eta_{\psi\circ \phi^{-1}}\|_{C^1}\le K\cdot 
\|\eta_\psi-\eta_\phi\|_{C^1}
$$
\end{lem}

\begin{proof} The Chain-rule for non-linearities 
$$
\eta_{\psi\circ \phi}(x)=\eta_{\psi}(\phi(x))\cdot D\phi(x) +\eta_\phi(x)
$$
implies
$$
\eta_{\phi^{-1}}(x)=- \eta_\phi(\phi^{-1}(x))\cdot D\phi^{-1}(x).
$$
Again the chain-rule gives
$$
\eta_{\psi\circ \phi^{-1}}=D\phi^{-1}  \cdot 
\left(\eta_\psi(\phi^{-1})-\eta_\phi(\phi^{-1}) \right). 
$$
Differentiation gives
$$
\begin{aligned}
D\eta_{\psi\circ \phi^{-1}}=&(D\phi^{-1})^2  \cdot 
\left(D\eta_\psi(\phi^{-1})-D\eta_\phi(\phi^{-1}) \right)\\
&+
D^2\phi^{-1}  \cdot 
\left(\eta_\psi(\phi^{-1})-\eta_\phi(\phi^{-1}) \right).
\end{aligned} 
$$
The bound  $\|\eta_\phi\|_{C^0}\le B$ gives a bound on  $\|\phi^{-1}\|_{C^2}$
and the Lemma follows.
\end{proof}

Now we are ready to prove the shuffling-Lemma \ref{shufflem}. 
The Lemmas \ref{Sandwich} and \ref{inv} imply the following estimate on the 
diffeomorphisms as defined in Lemma \ref{shufflem}
$$
\|\eta_{\Phi}-\eta_{\Phi^*}\|_{C^0}\le K 
\sum_{j=1}^n \|\eta_j-\eta^{*}_j\|_{C^0},
$$
where $K=K(B)$. One can integrate non-linearities and obtain
$$
\phi(x)=2\frac{\int_{-1}^x e^{\int_{-1}^s \eta_\phi} ds}
              {\int_{-1}^1 e^{\int_{-1}^s \eta_\phi} ds}-1.
$$
and
$$
D\phi(x)=2\frac{e^{\int_{-1}^x \eta_\phi} ds}
              {\int_{-1}^1 e^{\int_{-1}^s \eta_\phi} ds}.
$$
Notice that the Sandwich-Lemma \ref{Sandwich} implies that
$$
\|\eta_\Phi\|_{C^0}, \|\eta_{\Phi^*}\|_{C^0} \le K\cdot B.
$$
This uniform bound and the two expressions above can be used to get the 
desired 
estimate on the $C^2$ distance between $\Phi$ and $\Phi^*$ in \ref{shufflem}. 
We finished the proof of the Shuffling-Lemma.


\section{List of special notations}\label{list}

\begin{itemize}
\item [] $\beta_0$, $\beta_1$ saddle fixed points of a H\'enon-like map $F$, \S 3.4
\item [] $b=b_F$   average Jacobian of $F$, \S 6
\item [] $B_w^n=B_w^n(F)$  renormalization  pieces of level $n$, \S 5.2
\item [] $D^n_k$ derivative at the tip, \S 7.2
\item [] $F(x,y)= (f(x)-\eps(x,y), x)$ H\'enon-like map, \S 3.2 
\item [] $f_*$  fixed point of the unimodal renormalization $R_c$, \S 3.1
\item [] $f^*$  fixed point of the unimodal renormalization $R_v$, \S 3.1
\item [] $F_*$ fixed point of the  H\'enon-like renormalization, \S 4
\item[] $H$  non-linear part of coordinate change, \S 3.5
\item[] $\HH_\Omega$ space of analytic H\'enon-like maps, \S 3.3
\item[] $\II_\Omega(\overline{\epsilon})$ space of infintely renormalizable unimodal maps,
\item[] $\Jac F= |\di\eps/ \di y|$ Jacobian of $F$, \S 3.2  
\item[] $\lambda$ the universal scaling factor, \S 3.1
\item[] $\Lambda$ scaling part of coordinate change, \S 3.5
\item[] $\OO=\OO_F$ the critical Cantor set, \S 5.2
\item [] $\psi^1_v= H^{-1}\circ\La^{-1}$ coordinate change conjugating  $RF$ to $F^2$, \S 5.1
\item[] $\Psi_\omega^n=\Psi_\om^n(F)$ coordinate change conjugating  $R^n F$ to $F^{2^n}$, \S 5.1
\item[] $\Psi_k=\Psi_k^{k+1}=\Psi_v^1(R^k F)$, \S 7.2
\item[] $R_c$ renormalization operator near the ``critical point'', \S 3.1
\item[] $R_v$ renormalization operator near the ``critical value'', \S 3.1
\item[] $s_k$ non-linear part of the coordinate change $\Psi_k$, \S 7.2
\item[] $S^n_k$ non-linear part of the coordinate change $\Psi^n_k$, \S 7.2
\item[] $\sigma=\la^{-1}$ the universal scaling factor, \S 3.1
\item[] $t_k$ tilt, \S 7.2
\item[] $\tau=\tau_F$ tip, \S 7.2
\item[] $\UU_U$ space of analytic unimodal maps, \S 3.3
\item[] $v_*$ universal change of coordinates, \S 7.1

\end{itemize}

\end{document}